\documentstyle{article}

\newtheorem{theo+}           {Theorem}      [section]
\newtheorem{prop+}  [theo+]  {Proposition}
\newtheorem{coro+}  [theo+]  {Corollary}
\newtheorem{lemm+}  [theo+]  {Lemma}
\newtheorem{exam+}  [theo+]  {Example}
\newtheorem{rema+}  [theo+]  {Remark}
\newtheorem{defi+}  [theo+]  {Definition}

\newtheorem{exam+s}  [theo+]  {Examples}
\newtheorem{rema+s}  [theo+]  {Remarks}
\newtheorem{hyp+}  [theo+]  {Hypotheses}
\newtheorem{cla+}  [theo+]  {Claim}

\newenvironment{theorem}{\begin{theo+}}{\end{theo+}}
\newenvironment{proposition}{\begin{prop+}}{\end{prop+}}
\newenvironment{corollary}{\begin{coro+}}{\end{coro+}}
\newenvironment{lemma}{\begin{lemm+}}{\end{lemm+}}
\newenvironment{example}{\begin{exam+}\rm}{\end{exam+}}
\newenvironment{remark}{\begin{rema+}\rm}{\end{rema+}}
\newenvironment{definition}{\begin{defi+}\rm}{\end{defi+}}

\newenvironment{remarks}{\begin{rema+s}\rm}{\end{rema+s}}
\newenvironment{hypotheses}{\begin{hyp+}\rm}{\end{hyp+}}
\newenvironment{claim}{\begin{cla+}}{\end{cla+}}

\newcommand{\pa}{\partial}

\newcommand{\Bbb}{\bf}
\newcommand{\RR}{{\Bbb R}} 
\newcommand{\CC}{{\Bbb C}}  
\newcommand{\NN}{{\Bbb N}}  
\newcommand{\Rmm}{{\Bbb R}^{2m}} 
\newcommand{\Cmm}{{\Bbb C}^{2m}} 
\newcommand{\Cm}{{\Bbb C}^m} 
\newcommand{\Ck}{{\Bbb C}^k} 
\newcommand{\CP}{{{\Bbb C}P}^1}  

\newcommand{\JJ}{\cal J} 
\newcommand{\cc}{\Bbb c} 
\newcommand{\ii}{{\rm i}}  
 \newcommand{\ra}{\rightarrow}

\newcommand{\proof}{\noindent\ {\bf Proof} \hskip 0.4em}
\newcommand{\eproof}{\bigskip}

\begin{document}
\title{Hermitian structures and harmonic morphisms in higher dimensional
 Euclidean spaces}

\author{P. Baird and J. C. Wood
\thanks{Partially supported by EC grant CHRX-CT92-0050} \\ {\small
D\'epartement de Math\'ematiques, Universit\'e de Bretagne Occidentale} \\
{\small 6 Avenue Le Gorgeu, B.P. 452, 29275 Brest Cedex, France, and} \\
{\small Department of Pure Mathematics, University of Leeds}\\ {\small
Leeds LS2 9JT, G.B.} }

\date{}
\maketitle

\begin{abstract} We construct new complex-valued harmonic morphisms from
Euclidean spaces from functions which are holomorphic with respect to
Hermitian structures.  In particular, we give the first global examples
of complex-valued harmonic morphisms from $\RR^n$ for each
$n>4$ which do not arise from a K\"ahler structure; it is known that such
examples do not exist for $n \leq 4$. \end{abstract}

\section{Introduction}\label{sec:Introduction} A {\em harmonic map\/}
$\phi:U \ra \CC$ from an open subset $U$ of Euclidean $n$-space
$\RR^n$ is a solution of Laplace's equation
\begin{equation}
\Delta\phi \equiv \sum_{i=1}^n \frac{\pa^2 \phi}{(\pa x^i)^2} = 0, \quad
(x^1,\ldots x^n) \in U ;
\label{eq:Lap-Eucl}
\end{equation} such a map is a {\em harmonic morphism\/} if it
additionally satisfies the condition of {\em horizontal weak
conformality\/}:
\begin{equation}
\sum_{i=1}^n \left(\frac{\pa \phi}{\pa x^i}\right)^2 = 0, \quad
(x^1,\ldots x^n) \in U  . \label{eq:HWC}
\end{equation} (Equivalently, and in more generality, a harmonic morphism
is a map between Riemannian manifolds which pulls back germs of harmonic
functions to germs of harmonic functions --- see \cite{Fu-1,Is}.)

Suppose that $n$ is even, say $n = 2m$.  Then it is well-known that a map
$\phi:U \ra \CC$ which is holomorphic with respect to the standard
complex structure
$J$ on
$\Rmm$ (given, on identifying the latter with $\Cm$, by multiplication by
$\ii$) is a harmonic morphism; indeed $\phi$ is holomorphic if and only
if it satisfies the Cauchy-Riemann equations
$$
\frac{\pa\phi}{\pa x^{2j}} = \ii \frac{\pa\phi}{\pa x^{2j-1}}, \quad (j =
1,
\ldots m)
$$ which are easily seen to imply (\ref{eq:Lap-Eucl}) and (\ref{eq:HWC}).
More generally, if $\phi$ is holomorphic with respect to any K\"ahler
structure (i.e. orthogonal complex structure) on $U$, it is a harmonic
morphism --- this is an example of the general result that holomorphic
maps from a K\"ahler manifold to
$\CC$ (or to a Riemann surface) are harmonic morphisms \cite{Fu-1}.

Let us now give $U$ a non-K\"ahler (integrable) Hermitian structure (see
\S
\ref{sec:Hermitian structures}).  Then, unless
$J$ is cosymplectic (i.e. has zero divergence, see equation
(\ref{eq:divJ}) below), a holomorphic map $(U,J) \ra \CC$ still satisfies
(\ref{eq:HWC}), however it is no longer automatically harmonic. (Note
that it is, however, always a Hermitian harmonic map in the sense of Jost
and Yau \cite{Jost-Yau}.) In fact it is harmonic if and only if the
divergence of $J$ lies in the kernel of the differential of $\phi$  (see,
for example, \cite{Gu-Wo-2}).  It is not easy to see how to find maps
satisfying this condition in general.  However, in \cite{Wo-3} with some
further development in \cite{Ba-6}, strong relationships were shown
between Hermitian structures and harmonic maps from
$\RR^4$, (or, more generally, any $4$-dimensional anti-self-dual Einstein
manifold) to $\CC$ or to a Riemann surface; in particular, this provided a
way of finding all (submersive) harmonic morphisms from open subsets of
$\RR^4$ to
$\CC$.  In this paper we explore this relationship in higher dimensions,
and in particular give a new way of finding many interesting locally and
globally defined harmonic morphisms from open subsets of $\Rmm$ which are
holomorphic with respect to non-K\"ahler structures.

\medskip
 We find many interesting differences between the $4$-dimensional case
($m=2$) and the higher dimensional case ($m > 2$).  Our main results are
as follows:

(i)\quad In \cite{Wo-3} it is shown that in the $4$-dimensional case any
harmonic morphism
$\phi: \RR^4 \supset U
\ra
\CC$ is holomorphic with respect to precisely one Hermitian structure ---
unless it has totally geodesic fibres in which case it is holomorphic
with respect to precisely two Hermitian structures of opposite
orientations.  Further, in either case, the fibres are superminimal with
respect to these Hermitian structures, i.e. the Hermitian structures are
parallel (constant) along the fibres of $\phi$.  In the higher
dimensional case we construct harmonic morphisms ${\Rmm}\supset U\ra \CC$
which are holomorphic with respect to Hermitian structures, some with
fibres which are superminimal and others with fibres which are not
superminimal with respect to any Hermitian structure. The examples can be
chosen with
$J$ cosymplectic or not, see Theorem \ref{th:local-even}.  None of these
examples is holomorphic with respect to any K\"ahler structure and they
are {\em full}, (see Definition \ref{def:full}).    We give another
example which is holomorphic with respect to a family of Hermitian
structures and which has superminimal fibres with respect to some of the
family but not all.

(ii)\quad In the
$4$-dimensional case it is shown in
\cite{Wo-3} that a globally defined harmonic morphism $\RR^4 \ra \CC$ is
holomorphic with respect to a K\"ahler structure at least if it is
submersive.  In contrast, in the higher dimensional case we construct
globally defined harmonic morphisms which are not holomorphic with
respect to any K\"ahler structure, and in fact have fibres which are not
superminimal with respect to any almost complex structure, see Theorem
\ref{th:global-even}.  These can be chosen to be full or to factor to
full globally defined harmonic morphisms
$\RR^{2m-1} \ra \CC$, see Theorem \ref{th:global-odd}. They can be chosen
to be submersive or not.

(iii)\quad Thus, for all $n > 4$, we have full globally defined harmonic
morphisms
$\RR^n
\ra
\CC$ which do not arise (see Definition \ref{def:arise}) from K\"ahler
structures (Corollary \ref{cor:global-all}).

Since, by \cite {Ba-Ee}, the fibres of a submersive harmonic morphism to
$\CC$ (or to a Riemann surface) form a conformal foliation by minimal
submanifolds of codimension
$2$, we obtain many interesting such foliations, in particular (Theorem
\ref{cor:conf-all}), global foliations of $\RR^n$ for any  $n > 4$ which
do not arise from a K\"ahler structure.

Note that, for simplicity, our codomain is taken to be $\CC$.  However,
the methods of this paper apply replacing this by any Riemann surface. In
particular, in the last section, we discuss those harmonic morphisms
which factor to give harmonic morphisms from a $2n$-torus to a $2$-torus
and show, that, if $n>2$, as well as being holomophic with respect to
infinitely many K\"ahler structures, any such map is also holomorphic
with respect to infinitely many non-K\"ahler Hermitian structures.

In a sequel to this paper, we shall study which harmonic morphisms
constructed by our method factor to spheres, real, complex and
quaternionic projective spaces, giving new examples on these spaces. We
shall also consider harmonic morphisms which are holomorphic with respect
to more general CR structures on odd dimensional spaces.

\medskip We now describe the idea of our construction of harmonic
morphisms.

(1) Suppose that $\phi:\Rmm \supset U \ra \CC^k$ is holomorphic with
respect to an almost Hermitian structure $J$. Then its Laplacian, normally
a second order partial differential operator (\ref{eq:Lap-Eucl}) can be
expressed just in terms of first derivatives of $\phi$ and the matrix
valued function $M:U \ra \CC$ representing $J$, see Equation
(\ref{eq:Lap1}).

(2) Now suppose that $k=m$ and $z:\Rmm \supset U \ra \CC^k$ is
holomorphic with respect to a {\em Hermitian} (i.e. integrable almost
Hermitian) structure $J$ and has independent components. Then these
components give local complex coordinates for $(U,J)$ and so $M$ is a
function of them, holomorphic by integrability of $J$.  More generally,
this remains true for any $k$ and any holomorphic map $z$ if we suppose
that $J$ is constant along the components of the fibres of $z$. (In the
case
$k=1$ this means that $z:U \ra \CC$ has superminimal fibres.) In this
case the expression for the Laplacian of $z$ becomes (\ref{eq:Lap2}).

(3) Now, any holomorphic $z$ must be related to the standard local
complex coordinates $w = q-M(z)\bar{q}$ for $(U,J)$ by an equation of the
form (\ref{eq:F})  where $f$ is holomorphic.  In the case when $k=m$ we
can generally write $f$ in the simpler form (\ref{eq:f-can}).  In either
case, we can find a formula for the Laplacian of $z$ just in terms of the
holmorphic data $(f,M)$ (or $(h,M)$). Our idea is to choose this
holomorphic data such that one or more components $z^a$ of $z$ has zero
Laplacian and so is a harmonic mrophism. Examples \ref{ex:LSnotC} --
\ref{ex:family} are found this way.

(4) If two or more components $z^{a_1}, \ldots, z^{a_\ell}$ of $z$ are
harmonic, we may compose the function $(z^{a_1}, \ldots, z^{a_\ell}):U \ra
V \subset \CC^l$ with a holomorphic map $\psi:V
\ra \CC$ to obtain (see \S \ref{subsec:harhol}) more harmonic morphisms;
this device is used to find our global Examples \ref{ex:glob1},
\ref{ex:glob2}.

(5) We check that all examples $\psi\circ z$ are full and not holomorphic
with respect to a K\"ahler structure using tests in \S
\ref{subsec:fullsuperk}.

(6) All harmonic morphisms $U \ra \CC$ which are holomorphic with respect
to a Hermitian structure are constructed, at least locally, by our method,
see Proposition \ref{prop:main}.

\bigskip {\bf Acknowledgments}

Both authors would like to thank J. Eells and S. Gudmundsson for valuable
suggestions related to this work. The second author would like to thank
the D\'epartement de Math\'ematiques at the Universit\'e de Bretagne
Occidentale, Brest, for making possible a visit in July 1994 which
allowed this work to be completed.

\section{Parametrization of Hermitian structures on  $\RR^{2m}$.}
\label{sec:Hermitian structures}

Throughout this paper we shall adopt the double summation convention.

We shall consider the Euclidean space $\Rmm$ together with its standard
orientation.  Let $x\in \Rmm$. An {\em almost Hermitian structure at\/}
$x$ or {\em on\/} $T_x\Rmm$ is an orthogonal transformation $J_x:T_x\Rmm
\ra T_x\Rmm$, with the property that
$J_x^2 = -I$ where $I$ is the identity transformation.  We say that
$J_x$ is {\em compatible with the orientation\/} or {\em positive\/} if
there exists an oriented orthonormal basis of $T_x\Rmm$ of the form
$\{e_1, e_2, \ldots , e_{2m}\}$ with $e_{2j}=J_xe_{2j-1}, (j=1, \ldots
,m)$.  Otherwise, we say that $J_x$ is {\em incompatible with the
orientation\/} or {\em negative\/}.  By an {\em almost Hermitian
structure\/} $J$ on an open subset $U$ of $\Rmm$ we mean the choice of an
almost Hermitian structure at each point $x$ of $U$ which varies smoothly
with $x$.

By a {\em complex chart\/} $\phi :(V,J) \ra \Cm$ on an open  subset
$V$ of $U$ we mean a diffeomorphism of $V$ onto an open subset of
$\Cm$ such that $d\phi (JX) = \ii\, d\phi (X)$ for all $X\in T_x\Rmm$,
$x\in V$.  We say that $J$ is {\em integrable on\/} $U$ if we can find
such complex charts in the neighbourhood of each point of $U$.  Note that
this gives $U$ the structure of an
$m$-dimensional complex manifold.  An integrable almost Hermitian
structure on $U$ is also called a {\em Hermitian structure on \/} $U$. It
is called a {\em K\"ahler structure\/} or an {\em orthogonal complex
structure\/} if it is parallel, i.e. does not vary from point to point.
  The standard Hermitian (in fact K\"ahler) structure on $\Rmm$ is that
given by the complex chart $\Rmm\ra\Cm,\ (x^1,\ldots,x^{2m})\mapsto
(x^1+\ii x^2,\ldots,x^{2m-1}+\ii x^{2m})$.

Let $x\in\Rmm$ and suppose that $J_x$ is a positive almost Hermitian
structure at $x$.  Consider the complex linear extension
$J_x^{\cc}:T_x^{\cc}\Rmm\ra T_x^{\cc}\Rmm$ to the complexified tangent
space.  Since $(J_x^{\cc})^2 = -I$, the space
$T_x^{\cc}\Rmm$ splits into the sum of two eigenspaces
$T_x^{(1,0)}\Rmm$ and $T_x^{(0,1)}\Rmm$ corresponding to the eigenvalues
$+\ii$ and $-\ii$ of $J_x$ respectively; these are called the $(1,0)$-
and $(0,1)$-tangent spaces of $J_x$.  Then
$T_x^{(1,0)}\Rmm$ is spanned by vectors of the form $e-\ii J_x e\ (e\in
T_x\Rmm)$.  Such a vector  is {\em isotropic\/} in the sense that
$\langle e-\ii J_x e,\ e-\ii J_x e\rangle^{\cc} = 0$.  Here
$\langle \ , \ \rangle^{\cc}$ denotes the standard symmetric bilinear
inner product on $\CC^n$ given by $\langle z ,  w \rangle^{\cc} = z^1 w^1
+ \ldots  + z^n w^n$.

A complex subspace $V \subset \Cmm$ is called  {\em isotropic\/} if
$\langle v,v\rangle^{\cc} = 0$ for each $v\in V$.  (This implies that
$\langle v,w\rangle^{\cc} = 0$ for all $v,w\in V$.)  Such a subspace is
called {\em maximal (isotropic)} (or {\em Lagrangian\/}) if it is not
contained in a larger isotropic subspace; note that an isotropic subspace
is maximal if and only if it has complex dimension $m$.  Suppose that $V$
is maximal isotropic.  Then call
$V$ {\em positive\/} if it has a basis of the form $\{e_1-\ii e_2,\ldots
,e_{2m-1}-\ii e_{2m}\}$ where $\{e_1,\ldots ,e_{2m}\}$ is a positive
orthonormal basis of $\Rmm$.

For any $x\in \Rmm$ we may identify in a canonical way $T_x \Rmm$ with
$\Rmm$ and $T_x ^{\cc}\Rmm$ with $\Cmm$.  Then there is a one-to-one
correspondence between the set
${\JJ}_x^{+}(\Rmm)$ of positive almost Hermitian structures at $x$ and the
set
$\hbox{Iso}^{+}(\Cm)$ of positive maximal isotropic subspaces of
$T_x^{\cc}\Rmm$  given by sending $J_x\in {\JJ}_x^{+}(\Rmm)$ to its
$(1,0)$-tangent space.  Further, both spaces can be identified with the
Hermitian symmetric space $SO(2m)/U(m)$.  Indeed let $A \in SO(2m)$.  The
columns $\{ e_1,\ldots e_{2m} \}$ of $A$ define a positive orthonormal
basis.  Then the coset of $SO(2m)/U(m)$ defined by $A$ corresponds to the
positive almost Hermitian structure $J_x$ at $x$ with $J_x (e_{2j-1}) =
e_{2j}\ (j=1,\ldots m)$; this in turn corresponds to the maximal
isotropic subspace of
$\Cmm$ with basis  $\{e_1-\ii e_2,\ldots ,e_{2m-1}-\ii e_{2m}\}$.

We next define the {\em twistor bundle\/} $\pi :Z^{+}\ra\Rmm$  as the
bundle whose fibre at each point $x\in\Rmm$ consists of all positive
Hermitian structures at $x$.  As a real fibre bundle,
$Z^{+} = \Rmm \times SO(2m)/U(m)$.  The {\em vertical tangent bundle\/}
$T^V Z^{+}$ is the subbundle of $TZ^{+}$ consisting of tangents to the
fibres of $Z^{+}$, equivalently, $T^V Z^{+}$ is the kernel of the
differential $d\pi :TZ^{+}\ra T\Rmm$; the {\em horizontal tangent
bundle\/} $T^H Z^{+}$ is its orthogonal complement with respect to the
product metric on $Z^{+}$.   We give $Z^{+}$ a complex structure in a
standard way as follows:  For any $J\in Z^{+}$, we give
$T_J^V Z^{+}$ the almost Hermitian structure induced from that of
$SO(2m)/U(m)$.  For $T_J^H Z^{+}$, note that each $J\in Z^{+}$ defines an
almost complex structure on $T_{\pi(J)}\Rmm$; we pull this back to $T_J^H
Z$ via the isomorphism $d\pi_J$.  Finally, we define the almost complex
structure ${\JJ}$ on $Z^{+}$ by ${\JJ} = {\JJ}^H\oplus {\JJ}^V$.  As is
well-known, this structure is integrable.  Next note that an almost
complex structure $J$ on an open subset $U$ of
$R^{2m}$ can be considered as a map $J:U\ra SO(2m)/U(m)$ or as a section
$\sigma_J :U\ra Z^{+}$; then we have:

\begin{lemma}\cite{BdeB,Besse-Einstein}\label{le:int} Let $J$ be a
positive almost Hermitian structure on  $U\stackrel{open}{\subset}\Rmm$.
Then the following are equivalent:
\begin{enumerate}
\item $J$ is integrable;
\item the corresponding map $J:U\ra SO(2m)/U(m)$ is  holomorphic; \item
the corresponding section $\sigma_J:U\ra Z^{+}$ is  holomorphic;
\item the image of the section $\sigma_J$ is a complex  submanifold of
$(Z^{+},{\JJ})$.
\end{enumerate}
\end{lemma}

We must now interpret our development in terms of the  {\em co\/}-tangent
space.  Note that the natural linear isomorphism
$\bar{\flat} :T_x^{\cc}\Rmm \ra T_x^{\star \cc}\Rmm$ given by
$X^{\bar{\flat}}(Y) =
\langle \bar{X},  Y\rangle^{\cc},\ (X,Y\in T_x^{\cc}\Rmm)$ gives a
one-to-one correspondence between (positive) isotropic subspaces of
$T_x^{\cc}\Rmm$ and those of  $T_x^{\star \cc}\Rmm$; the subspace of
$T_x^{\star \cc}\Rmm$  corresponding to the $(1,0)$-tangent space (resp.
$(0,1)$-tangent space) of a positive almost Hermitian structure  being
called its {\em $(1,0)$-cotangent space\/},
$T_x^{\star (1,0)}\Rmm$ (resp. its {\em $(0,1)$-cotangent space\/},
$T_x^{\star (0,1)}\Rmm$).

Next, suppose that $V_0\in \hbox{Iso}^{+}(\Cmm)$ is the
$(1,0)$-cotangent space of a positive almost Hermitian structure $J_0\in
{\JJ}_x^{+}(\Rmm)$.  Note that $\Cmm = V_0\oplus\bar{V_0}$.  Then it is
easily seen that the graph of a linear map $M:V_0\ra\bar{V_0}$ belongs to
$\hbox{Iso}^{+}(\Cmm)$ if and only if $M$ is skew-symmetric.  This
correspondence provides a chart for
$\hbox{Iso}^{+}(\Cmm)$ centred on $V_0$ or, equivalently, for
${\JJ}_x^{+}(\Rmm)$ centred on $J_0$.  Explicitly, let $\{q^j =
x^{2j-1}+\ii x^{2j}:j=1\ldots m\}$ denote complex coordinates for
$(\Rmm,J_0)$, then the $(1,0)$- and $(0,1)$- cotangent spaces of
$J_0$ have bases $\{dq^j: j=1,\ldots m\}$ and $\{dq^{\bar{j}}: j=1,\ldots
m\}$ respectively.  The above chart centred on $J_0$ can be described as
follows:

Define a map
\begin{equation}\label{eq:chart}
\iota: \CC^{m(m-1)/2}  \stackrel{\cong}{\ra} so(m,\CC)\ra
Iso^{+}(\Cmm)\ra {\JJ}_x^{+}(\Rmm) \label{eq:iota}
\end{equation} by
\begin{equation}
\mu = (\mu_1,\ldots ,\mu_{m(m-1)/2}) \mapsto M \mapsto V \mapsto J
\label{eq:param}
\end{equation}
 where  $M = M(\mu)$ is the $m\times m$ antisymmetric matrix given by
$$
\big(M^i_{{}{\bar{j}}}(\mu)\big) =
\left(\begin{array}{rrrrr} 0          & {}\mu_1     & {}\mu_2    &
\cdots & \mu_{m-1} \\ -\mu_1    & 0           & {}\mu_m    & \cdots &
\mu_{2m-3} \\ -\mu_2     & -\mu_m      & 0          & \cdots &
\mu_{3m-6} \\ \vdots     & \vdots      & \vdots     & \ddots &
\vdots   \\ -\mu_{m-1} & -\mu_{2m-3} & -\mu_{3m-6}  & \cdots & 0
\end{array}\right) \label{M},
$$
$V = V(M)$ is the positive maximal isotropic subspace given by
$$ V = \hbox{span}\{ e^{i} = dq^{i}-M^{i}_{{}{\bar{j}}} dq^{\bar{j}}  :
i=1,\ldots,m\}
$$ and $J = J(V)$ is the positive Hermitian structure at $x$ with
$(1,0)$-cotangent space
$V$. Note that these mappings are well-defined since:

(i) we have
\begin{eqnarray*}
 \langle e^{i},e^{j} \rangle^{\cc}  & = &
\langle  dq^{i} - M^{i}_{{}{\bar{k}}}dq^{\bar{k}} , dq^{j} -
M^{j}_{{}{\bar{\ell}}}dq^{\bar{\ell}} \rangle^{\cc} \\ & = &
2(M^{i}_{{}{\bar{j}}} + M^{j}_{{}{\bar{i}}}) \\ & = & 0
\end{eqnarray*} since $M$ is skew-symmetric and so $V = V(M)$ is
isotropic, further the $e^{i}$ form a positively oriented basis and so
$V$ is maximal isotropic and positive;

(ii) any positive maximal isotropic subspace $V$ is the $(1,0)$-cotangent
space of a unique positive almost Hermitian structure
$J = J(V)$ at the point $(q^1,\ldots,q^m)\in \Cm = \Rmm$.

Note that the map (\ref{eq:chart}) gives a chart for a dense subset of
$\hbox{Iso}^{+}(\Cmm) = SO(2m)/U(m)$  called {\em a large cell}  --- in
fact it is essentially the exponential map.

For simplicity of notation we shall frequently indicate the various
compositions of maps in (\ref{eq:param}) by giving the dependent and
independent variables, e.g.  $J(M) = J(V(M))$ and $J(\mu) =
J(V(M(\mu)))$, (i.e.
$\iota$), and their inverse maps by $M(J)$, $\mu (J)$ etc.

For computations, the following is useful (cf. \cite{New-Nir}):

\begin{lemma}  The $(0,1)$-tangent space of $J_x$ has basis
$\{e_{\bar{i}}:i=1,\ldots,m\}$ where
$$ e_{\bar{i}} = \frac{\pa}{\pa q^{\bar{i}}} +
M^{j}_{\bar{i}}\frac{\pa}{\pa q^{j}}.
$$
\end{lemma}
\proof It suffices to note that, under the natural pairing of
$T_x^{\star\cc}(\Rmm)$ and $T_x^{\cc}(\Rmm)$ the vectors
${e^{i}}$ annihilate the vectors $e_{\bar{j}}$.
\eproof

\begin{example}\label{ex:M-m=2} Let $m=2$.  Then
$$ M(\mu) = \left( \begin{array}{rr} 0   &   \mu\\  -\mu & 0
\end{array}\right)
$$ and the almost Hermitian structure $J(\mu)$  corresponding  to
$\mu\in\CC$ has $(1,0)$-cotangent space spanned by $\{dq^{1}-\mu
dq^{\bar{2}},\ dq^{2}+\mu dq^{\bar{1}}\}$.  This extends to a
diffeomorphism $\CP\ra \hbox{Iso}^{+}(\RR^4) \cong SO(4)/U(2)$ given by
$$ [\mu_0,\mu_1] \mapsto \hbox{span}\{\mu_0 dq^1-\mu_1 dq^{\bar{2}},
\mu_0 dq^2+\mu_1 dq^{\bar{1}}\}.
$$
\end{example}

\begin{example}  Let $m=3$. Then, attempting to extend the map
(\ref{eq:param}) to $\CC P^3$ by introducing $\mu_0$ as in Example
\ref{ex:M-m=2} we find that the corresponding covectors
$\mu_0 dq^1 - \mu_1 dq^{\bar{2}} - \mu_2 dq^{\bar{3}},\ \mu_0 dq^2 + \mu_1
dq^{\bar{1}} - \mu_3 dq^{\bar{3}},\ \mu_0 dq^3 + \mu_2 dq^{\bar{1}} +
\mu_3 dq^{\bar{2}}$ are linearly dependent if $\mu_0 = 0$.  However,
consider the fourth covector: $\mu_3 dq^1 - \mu_2 dq^2 +
\mu_1 dq^1$. If $\mu_0 \neq 0$, this is a linear combination of the above
covectors, if $\mu_0 = 0$, it isn't; in both cases the {\em four\/}
covectors span a positive isotropic subspace of $\CC^6$ of dimension {\em
three\/} and so define a diffeomorphism  $\CC P^3 \ra
\hbox{Iso}^{+}(\CC^6) \cong SO(6)/U(3)$ extending (\ref{eq:param}).
\end{example}

Note that in higher dimensions, no such extension of (\ref{eq:param}) is
possible since then  $\hbox{Iso}^{+}(\Rmm)$ is not homeomorphic to any
complex projective space.

\bigskip We next construct complex coordinates on the twistor bundle
$Z^{+} = \Cm\times{\JJ}^{+}(\Rmm)$.  Note that this has {\em product
coordinates\/} given by the formula $(q,J) \mapsto (q^1,\ldots,
q^m,\mu_1(J) ,\ldots, \mu_{m(m-1)/2}(J))$ but these are not complex
coordinates for the complex structure ${\JJ}$ on $Z^{+}$ since the $q^i$
are not holomorphic functions on $Z^{+}$.  However we have

\begin{proposition} The functions
$(q,\mu)\mapsto w^i = q^i - M^i_{{}{\bar{j}}}(\mu)q^{\bar{j}},\
(i=1,\ldots, m)$ and $(q, \mu)\mapsto \mu^i, \ (i=1,\ldots, m(m-1)/2)$
define complex coordinates on a dense open subset of the twistor space
$Z^{+} = \Rmm\times{\JJ}^{+}(\Rmm)$.
\end{proposition}

\proof   It is clear that the $\mu^i$ are holomorphic on
$Z^{+}$ (indeed they give complex coordinates for each fibre), whilst
\begin{eqnarray*} dw^i & = & d(q^i - M^i_{{}{\bar{j}}}(\mu)q^{\bar{j}})
\\ & = & dq^i - M^i_{{\bar{j}}}(\mu)dq^{\bar{j}} +
q^{\bar{j}}dM^i_{{}{\bar{j}}}(\mu)
\end{eqnarray*} is the sum of covectors of type $(1,0)$, since
$M^i_{{\bar{j}}}(\mu)$ is holomorphic in $\mu$.  Thus the $w^i$ are also
holomorphic.  (Alternatively we calculate
$e_{\bar{j}}(w^i) = 0$.) Clearly the set of differentials
$\{dw^i:i=1,\dots m\}\cup \{dq^j:j=1,\dots m(m-1)/2\}$ is linearly
dependent on a dense set
\eproof

We can also interpret the $w^i$ as complex coordinates  on
$\Rmm$ as follows:

\begin{corollary}\label{cor:w-cxcoords} Let $J = J(M(\mu(q)))$ be a
Hermitian structure on an open set
$U$ of $\Rmm$. Then the functions $w^i:U \ra \CC$, $q\mapsto w^i = q^i -
M^i_{{\bar{j}}}(\mu(q))q^{\bar{j}}$ are complex coordinates for $(U,J)$.
\end{corollary}

\proof  The functions  $q\mapsto w^i$ are compositions of the
holomorphic coordinate functions $(q,\mu)\mapsto w^i$ on $Z^+$ and  the
holomorphic section $\sigma_J$ corresponding to $J$ (see Lemma
\ref{le:int}).
\eproof

We now turn to the construction of maps which are holomorphic with
respect to a Hermitian structure.  In \cite{Wo-3} harmonic morphisms
$\phi$ from open sets $U$ of $\RR^4$ to $\CC$ (or a Riemann surface) are
constructed which are holomorphic with respect to a Hermitian structure
$J$ which is constant along the fibres of
$\phi$ and is thus given by a composition of the form
$J:U\stackrel{\phi}\ra V
\stackrel{\mu}{\ra}\CC$ where $V$ is open in $\CC$ and $\mu$ is
holomorphic.   In the next section we generalise this to consider maps
$\phi:\Rmm\supset U\ra\Ck$ which are holomorphic with respect to a
Hermitian structure $J$ that is constant along the fibres of $\phi$ and
thus is given by the composition of $\phi$ with a holomorphic map

\begin{equation}
\CC^k \stackrel{open}\supset V  \ra so(m,\CC),\quad z \mapsto M(z).
\label{eq:M}
\end{equation}

Suppose that $\phi:\Rmm \stackrel{open}{\supset} U \ra V \subset \Ck$,\
$z=\phi (q)$, is given implicitly in terms of the complex coordinates $w$
for $(U,J)$ by an equation
\begin{equation} f(w,z) = 0 \label{eq:f}
\end{equation} where
$f:\Cm\times\Ck \stackrel{open}{\supset} W \ra \Ck $ is holomorphic.
Then, in terms of standard coordinates
$q=(q^1,\ldots,q^m)$ on $\Cm=\Rmm$ , the map $z=\phi(q)$ is a solution to
an equation of the form

\begin{equation} f(q-M(z)\bar q, z) = 0 \label{eq:f2}
\end{equation} where $M: \CC \stackrel{open}\supset V
\ra so(m,\CC) = \CC^{m(m-1)/2}$ is a given holomorphic function.  In the
next section we consider when the components of a solution of such an
equation are harmonic and thus harmonic morphisms.

\section{Harmonic holomorphic mappings}
\label{sec:Haho_maps} In this section,  we calculate the Laplacian of a
map from an open subset $U$ of $\Rmm$  which is holomorphic with respect
to a Hermitian structure on $U$, firstly (\S
\ref{subsec:Laphol}) in general, then specializing first to Hermitian
structures constant along the fibres (\S \ref{subsec:Lapholsup}) then (\S
\ref{subsec:Lapholimpl})  to holomorphic maps defined implicitly by
Equation (\ref{eq:f2}). Then in Section
\ref{sec:examples} we use these formulae to find examples of harmonic
holomorphic maps.  On the way, we make some general remarks about the
space of all harmonic holomorphic maps on $U$ (\S \ref{subsec:harhol}).

\subsection{Formula for the Laplacian of a map which is holomorphic with
respect to some  Hermitian structure}\label{subsec:Laphol}

Let $U$ be an open subset of $\Rmm$ and let $J:U\ra{\JJ}^{+}({\Rmm})$, $J
= J(M(q))$ be a positive almost Hermitian structure.

Let $\phi:U\ra {\Ck}$ be a smooth map.  We wish to calculate the
Laplacian of $\phi$. Now $\phi$ is holomorphic with respect to $J$ if and
only if
\begin{equation}
\frac{\partial\phi}{\partial q^{\bar{i}}}+M^j_{{}{\bar{i}}}(q)
\frac{\partial\phi}{\partial q^j} = 0 \quad \hbox{for } i = 1,\ldots , m
. \label{eq:hol}
\end{equation} Assuming that this is the case,  differentiating
(\ref{eq:hol}) with respect to $q^i$ and summing over $i$ gives
$$
\frac{1}{4} \Delta\phi \equiv \sum_{i=1}^{m} \frac{\partial
^{2}\phi}{\partial  q^{i}\partial q^{\bar{i}}}  = -  \sum_{i,j=1}^{m}
M^{j}_{{}{\bar{i}}}
\frac{\partial^{2} \phi }{\partial q^{i} \partial q^{j}}  -
\sum_{i,j=1}^{m} \frac{\partial M^{j}_{{}{\bar{i}}}}{\partial q^{i}}
\frac{\partial \phi}{\partial q^{j}}.
$$

Now the first term is zero by the antisymmetry of $M^{j}_{\bar{i}}$ and
the symmetry of
$(\partial ^{2}\phi /\partial q^{i}\partial q^{j})$, hence we get

\begin{lemma}\label{Laplacian1} Let $\phi :U\ra {\Ck}$ be holomorphic
with respect to an almost Hermitian structure
$J = J(M(q))$ on $V$.  Then
\begin{equation}
\frac{1}{4} \Delta\phi = - \sum_{i,j=1}^{m} \frac{\partial
M^{i}_{\bar{j}}}{\partial q^{j}}
\frac{\partial \phi}{\partial q^{i}}. \label{eq:Lap1}
\end{equation}
\end{lemma}

\begin{remarks}\label{re:div}  1.\ We do not need $J$ to be integrable
for this lemma.

2.\ More generally, for any oriented Riemannian manifold $M^{2m}$ and map
$\phi :M^{2m}\ra {\Ck}$ holomorphic with respect to an almost Hermitian
structure $J$ on $M^m$ we have \cite{Gu-Wo-2}:

\begin{equation}
\Delta \phi = -d\phi (J\delta J). \label{eq:Lap-divJ}
\end{equation}

3.\  Formula (\ref {eq:Lap1}) breaks down where $J(q)$ does not belong to
a large cell for all $q \in V$ and so cannot be written as $J(M(q))$; but
unless
$J$ is constant, such points are nowhere dense and, if $J$ is integrable,
they are isolated since then $J:V\ra {\JJ}^{+}({\Rmm})$ is holomorphic.
We can thus replace $V$ by a slightly smaller open set $V'$ such that
$J(V')$ does belong to a large cell.
\end{remarks}

\subsection{The space of all harmonic holomorphic  maps and cosymplectic
manifolds}\label{subsec:harhol} Let $U$  be an open subset of $\Rmm$,\
($m \in \NN$) and let $J$ be an almost Hermitian structure on $U$.
\begin{proposition}\label{prop:sethaho} The set  $\hbox{\rm Haho}(U) = \{
\phi:U
\ra \CC^p: p \in \NN,\ \phi$ {\em is harmonic and holomorphic with
respect to} $J\}$ is closed under postcomposition with holomorphic maps
$f:\CC^p \ra \CC^q$, i.e. if $\phi
\in \hbox{\rm Haho}(U)$ and $f:V \ra \CC^q$ is holomorphic, where
$V$ is an open subset of $\CC^p$ containing $\phi(U)$ then the composition
$f \circ \phi$ is in $\hbox{\rm Haho}(U)$.
\end{proposition}
\proof  By the chain rule, the Laplacian of such a composition $f \circ
\phi$ is given by
\begin{eqnarray*}
\Delta (f \circ \phi)  & = & \frac{\pa^2f}{\pa z^j \pa z^k} \langle \nabla
\phi^j, \nabla \phi^k \rangle^{\cc} + \frac{\pa f}{\pa z^l} \Delta \phi^l
\\
                       & = & 0
\end{eqnarray*} where $\nabla \phi^j = (\pa \phi^j / \pa x^1, \ldots , \pa
\phi^j / \pa x^{2m})$ for $(x^1, \ldots , x^{2m}) \in \Rmm$.  Now, as in
\cite{Gu-4},
$\nabla \phi^j$, $\nabla \phi^k$ and $\nabla \phi^j + \nabla \phi^k$,
being the gradients of holomorphic functions, are isotropic, so by the
polarization identity, $\langle \nabla \phi^j, \nabla \phi^k
\rangle^{\cc} = 0$.
\eproof

There is one special case when {\em all\/} holomorphic maps $\phi:U \ra
\CC^p$ are harmonic, namely when $J$ is {\em cosymplectic}, i.e. its
divergence $\delta J$ is zero, $\delta J$ being defined by
\begin{equation}
\delta J=\sum (\pa_{e_{i}}J)(e_{i}) \label{eq:divJ}
\end{equation} where $\{e_i\}$ is any orthonormal basis of the tangent
space (see \cite{Lichn, Gu-Wo-2}).

Precisely we have
\begin{proposition}\label{prop:cosymplectic}
 Let $U$ be an open subset of $\Rmm$ and let
$J = J(M)$ be a Hermitian structure on it.  Then the following are
equivalent:

(1) for any point $x \in U$ there are $m$ holomorphic harmonic maps $U'
\ra \CC$ defined on a neighbourhood $U'$ of
$p$ which are {\em independent}, i.e. their gradients are linearly
independent over $\CC$;

(2) all holomorphic maps $U' \ra \CC$ defined on an open subset $U'$ of
$U$ are harmonic;

(3) the complex coordinates
$w^i = q^i - M^i_{\bar{j}}(q) q^{\bar{j}}$ (see Corollary
\ref{cor:w-cxcoords}) are harmonic;

(4)  the matrix $M$ satisfies
$$
\sum_{j = 1}^m
\frac{\pa M^i_{\bar{j}}}{\pa q^j} = 0 \quad \hbox {for all}\quad i = 1
\ldots m;
\label{eq:cosym}
$$

(5)  The Hermitian structure $J$ is cosymplectic.
\end{proposition}

\proof  That (5) implies (2) follows from the formula
(\ref{eq:Lap-divJ}), and is, in any case, well-known (see
\cite{Lichn}).  For the converse, apply this formula to $m$ independent
holomorphic harmonic functions to see that all the components of the
divergence of $J$ vanish.

That (1) implies (2) follows by noting that, if $\phi_1, \ldots ,\phi_m
:U' \ra \CC $ are independent holomorphic functions, then any holomorphic
function $\phi:U' \ra \CC$ is, locally a function of them, and so, by
Proposition \ref{prop:sethaho} , is harmonic.
 The converse is trivial and the equivalence of (2) and (3) is similar.
That (3) implies (4) follows from Formula (\ref{eq:Lap1}) applied to
$\phi = w^k$.  Similarly, that (4) implies (5) follows from Formula
(\ref{eq:Lap-divJ}) applied to $\phi = w^k$.

\subsection{The Laplacian of a holomorphic map when the Hermitian
structure is constant along the fibres}\label{subsec:Lapholsup}

Let $z:\Rmm \supset U \ra \CC,\ q \mapsto z(q)$ be a map which is
holomorphic with respect to a Hermitian structure $J = J(M(q))$ on U.
Suppose that $J$ is constant along the components of the fibres of $z(q)$.
(Note that this is automatically true if $k=m$ and $z(q)$ has independent
components. For then these components give local complex coordinates for
$(U,J)$ and so $M$ is a funtion of them.)    In this case $J$ is of the
form:

\begin{equation} J:{\Rmm}\supset U \ra V \ra so(m,\CC)\ra
{\JJ}^{+}({\Rmm}) \label{eq:J-const-fibres}\ ,\quad q \mapsto z \mapsto M
\end{equation} for some map $\Ck \stackrel{open}{\supset}V\ra so(m,
\CC)$, \ $z \mapsto M(z)$ which is holomorphic if
$J$ is integrable.  Assuming that is the case, by the chain rule  Formula
(\ref{eq:Lap1}) becomes

\begin{equation}
\frac{1}{4} \Delta z ^{a} = -\sum_{i,j=1}^{m}\sum_{b=1}^{k}
A^{a}_{{}i}A^{b}_{{}j}\frac{\partial M^{i}_{\bar{j}}}{\partial z^{b}}\
,\quad  (a=1,....,k)\ ,\label{eq:Lap2}
\end{equation} where
\begin{equation} A^{a}_{i}=\partial z^{a}/\partial q^{i}\ ,\quad
(a=1,...k,\ i=1,...m), \label{eq:A}
\end{equation} which we may write in the form

\begin{equation}
\frac{1}{4}\Delta z^{a} = -\sum_{1\leq i<j\leq m}^{}\sum_{b=1}^k
C^{ab}_{ij}
\frac{\partial M^{i}_{\bar{j}}}{\partial z^{b}}\ , \quad (a = 1, \ldots
k)
\label{eq:Lap2'}
\end{equation} where
\begin{equation} C^{ab}_{ij}=A^{a}_{i}A^{b}_{j}-A^{b}_{i}A^{a}_{j},\
(i,j=1,\ldots,m,\ ,\quad a,b=1,\ldots,k) \label{eq:C}.
\end{equation} (Note that $C^{ab}_{ij}=-C^{ab}_{ji} =-C^{ba}_{ij}$, in
particular $C^{aa}_{ij}=C^{ab}_{ii}=0$.)  We thus have proved the

\begin{lemma}\label{le:Lap-Jconst} Suppose that $z:U \ra \Ck$,\ $q
\mapsto z(q)$ is holomorphic with respect to a Hermitian structure $J$
given by a composition of the form (\ref{eq:J-const-fibres}), thus $J(q)
= J(M(z(q)))$ for $q \in U$.  Then the Laplacian of $z(q)$ is given by
equations (\ref{eq:Lap2}, \ref{eq:A}), or equivalently by equations
(\ref{eq:Lap2'}, \ref{eq:C}).
\end{lemma}

Our first example involves the concept of superminimality:
\begin{definition}\label{def:superminimal} Let $S$ be an even dimensional
submanifold of an open subset $U$ of $\Rmm$.  We say that an almost
complex structure $J$ on $U$ is {\em adapted\/} to $S$ if
$S$ is a complex submanifold with respect to $J$, i.e. $TS$ is closed
under $J$.  If
$S$ is adapted, we say that it is {\em superminimal with respect to
$J$} if $J$ is parallel (constant) along the submanifold.
\end{definition}

If $J$ is an almost Hermitian structure, superminimality implies
minimality, see for example, \cite{Gu-Wo-2}.

\begin{example} Suppose that $k=1$, that is, suppose that $z:\Rmm \supset
U \ra \CC$,\
$q \mapsto z(q)$ is holomorphic  with respect to a Hermitian structure
$J$ that is constant along the  fibres of $z(q)$; this implies that the
regular fibres of $z(q)$ are superminimal.  Then (\ref{eq:Lap2'}) reduces
to
$$
\Delta z =0.
$$
\end{example} This is clear since $C^{11}_{ij}=0$ is the only coefficient.
Of course we know this since superminimal fibres are minimal so that $z$,
being horizontally weakly conformal, is a harmonic morphism by
\cite{Ba-Ee}.

\begin{remark} It can be seen directly that, for a conformal $J$-closed
foliation of codimension $2$ of a Hermitian manifold $(M^{2m}, J)$,  the
condition $\nabla _{V}J=0$ for all $V$ tangent to the leaves implies that
 $\delta J$ is tangent to the leaves (see \cite{Gu-Wo-2}).
\end{remark}

\subsection{The Laplacian of an implicitly defined holomorphic
map}\label{subsec:Lapholimpl}

As at the end of \S 2, we consider holomorphic mappings $z :(V, J)\ra
{\Ck}$,\
$z=z(q)$ given implicitly by an equation of the form (\ref{eq:f2}).  To be
precise our hypotheses are as follows:

\begin{hypotheses}\label{hyp:general} We suppose that
$f:{\Cm}\times {\Ck}\stackrel{open}{\supset}W\ra {\Ck}$ and
$M :{\Ck}\stackrel{open}{\supset}V\ra so(m,\CC)$ are given holomorphic
mappings.

Defining $F:{\Rmm}\times {\Ck}\stackrel{open}{\supset}A\ra {\Ck}$  by
$F(q,z)=f(q-M(z)\bar{q},z)$ on the open set $A=\{(q,z): z\in V,\
(q-M(z)\bar{q}, z)\in W\}$, supposed non-empty, we suppose that the
Jacobian matrix
\begin{equation} K=(\partial F^{a}/\partial z^{b})_{a,b=1,\cdots
,m}\label{eq:K}
\end{equation} is non-singular  on a non-empty open subset
$A'\subset A$.

Lastly we suppose that $z :{\Rmm}
\stackrel{open}{\supset}U\ra {\Ck},\ z=z(q)$ is a smooth local solution
to the equation
\begin{equation} F(q,z) \equiv f(q-M(z)\bar{q},z)  = 0 \label{eq:F}
\end{equation} (i.e.\ $z(U)\subset V$,\ $(q,z(q))\in A'$ for all $q\in U$,
and
\begin{equation}  F(q, z(q))=0 \label{eq:F2}
\end{equation} for all $q \in U$.)
\end{hypotheses}

We can now give a formula for the Laplacian of such a solution to
equation (\ref{eq:F}):

\begin{proposition}\label{th:super} Suppose that Hypotheses
\ref{hyp:general} hold, in particular, suppose that
$z: U \ra \Ck , z=z(q)$ is a smooth local solution to equation
(\ref{eq:F}).  Then

(i)  $z :U\ra {\Ck}$,\ $q \mapsto z(q)$ is holomorphic with respect to
the Hermitian structure $J$ on $U$ given by $J(q)=J(M(z(q)))$ for $q\in
U$;

(ii)  this Hermitian structure is parallel (i.e. constant) along the
fibres of $z(q)$;

(iii)  the Laplacian of each component of $z(q)$ is given by
(\ref{eq:Lap2}),  or equivalently (\ref{eq:Lap2'}, \ref{eq:C}), where
\begin{equation} A^{a}_{i}= - (K^{-1})^{a}_{b} \frac{\partial
f^{b}}{\partial w^{i}}. \label{eq:A2}
\end{equation}
\end{proposition}

\begin{remarks} 1.\ That (\ref{eq:F}) has smooth local solutions $z(q)$
is guaranteed by the Implicit Function Theorem; indeed given any
 $(q_{0},z_{0})\in A$, (\ref{eq:F}) has a unique local smooth solution
$z:U\ra {\Ck},\ z=z(q)$, on some neighbourhood $U$ of $q_{0}$, such  that
$z(q_{0})=z_{0}$.

2.\ That $J$ is parallel along the fibres of $z(q)$ does not of course
imply that $J$ is parallel along the fibres of each component $z^{a}$ of
$z(q)$.
\end{remarks}

\proof  (i)  The map $z(q)$ is of the form
$$ {\Rmm}\supset U \ra {\Cm}\stackrel{\psi}{\ra}{\Ck},
\quad q\mapsto w\mapsto z
$$   where, using matrix notation for convenience, $w(q) =
q-M(z(q))\bar{q}$, and $z = \psi (w)$ is a local solution to $f(w,z)=0$,
i.e. $f(w, \psi (w))\equiv 0$.   By the Implicit Function Theorem, $\psi$
is holomorphic.  Now since, by Corollary
\ref{cor:w-cxcoords} the components of $w$ give local coordinates,
$w(q)$ is  holomorphic with respect to the almost Hermitian structure $J$
on $U$ given  by the composition
\begin{equation} J:U\stackrel{z}{\ra}V\subset
{\Ck}\stackrel{M}{\ra}so(m,\CC)
\cong C^{m(m-1)/2} \stackrel{\iota}{\ra} {\JJ}^{+}({\Rmm}).\label{eq:J}
\end{equation} Thus $z(q)$ is holomorphic with respect to $J$.

(Alternatively it can be seen from the chain rule for $F(q,z)=f(w(q,z),z)$
where $w(q,z)=q-M(z)\bar{q}$ that, for $\alpha = 1, \ldots k$,
$$
\frac{\partial F^{\alpha}}{\partial q^{i}}=\frac{\partial f^{\alpha}}
{\partial w^{i}} \quad\hbox{and}\quad \frac{\partial F^{\alpha}}{\partial
q^{\bar{i}}}= -M^{j}_{\bar{i}}(z)\frac{\partial f^{\alpha}}{\partial
w^{j}}
$$ so that
\begin{equation} (\frac{\partial}{\partial
q^{\bar{i}}}+M^{j}_{\bar{i}}(z)\frac{\partial } {\partial
q^{j}})F^{\alpha}=0, \label{eq:holF}
\end{equation} expressing the holomorphicity of $q\mapsto F(q,z)$ with
respect to $J$. Then, differentiating (\ref{eq:F}) with respect to
$q^{I},\ (I =1, \cdots ,\bar{m})$ gives
$$
\frac{\partial F^{\alpha}}{\partial q^{I}}+\frac{\partial F^{\alpha}}
{\partial z^{a}}\frac{\partial z^{a}}{\partial q^{I}} = 0
$$ so that
$$
\frac{\partial z^{a}}{\partial q^{I}}=-(K^{-1})^{a}_{\alpha}
\frac{\partial F^{\alpha}}{\partial q^{I}},
$$ whence from (\ref{eq:holF}) we obtain
$$ (\frac{\partial}{\partial
q^{\bar{i}}}+M^{j}_{\bar{i}}(z)\frac{\partial } {\partial q^{j}})z^{a}=0,
$$ showing that $z(q)$ is holomorphic with respect to $J$.)

That $J$ is Hermitian follows from the fact that (\ref{eq:J}) exhibits $J$
as the composition of holomorphic maps, so it is a holomorphic map
 $(U, J)\ra {\JJ}^{+}({\Rmm})$.  By Lemma \ref{le:int} it is integrable
and so  Hermitian.

(ii)  Obvious from the form of $J$ in (\ref{eq:J}) .

(iii)  By Lemma \ref{le:Lap-Jconst} the Laplacian of $z(q)$ is given by
equations (\ref{eq:Lap2}, \ref{eq:A}).  Now, differentiating
(\ref{eq:F2}) totally with respect to $q^i$, we get, for $i = 1, \ldots ,
m$,
\begin{equation}
\frac{\partial F^{b}}{\partial q^{i}}+\frac{\partial F^{b}}{\partial
z^{a}}
\frac{\partial z^{a}}{\partial q^{i}}=0 \label{eq:F-q}
\end{equation} whence
\begin{equation}
\frac{\partial z^{a}}{\partial q^{i}}=-(K^{-1})^{a}_{b}\frac{
\partial F^{b}}{\partial q^{i}}. \label{eq:dzdq}
\end{equation}
$F(q,z)=f(q-M(z)\bar{q},z)$,
$$
\frac{\partial F^{b}}{\partial q^{i}}=\frac{\partial f^{b}}{\partial
w^{a}}
\frac{\partial w^{a}}{\partial q^{i}}=\frac{\partial f^{b}}{\partial
w^{i}}
$$ giving the formula (\ref{eq:A2}) for the matrix $A$.
\eproof

\bigskip {\em The special case $k=1$}.

In this case the right hand side of (\ref{eq:Lap2'})  is always zero and
we get a result:

\begin{proposition} Suppose that Hypotheses \ref{hyp:general} hold with
$k = 1$, in particular, suppose that
$z:\Rmm \supset U \ra \Ck$,\ $q \mapsto z(q)$ is a smooth local solution
to equation (\ref{eq:F}) with
$k=1$.  Then

(i)  $z:U\ra {\CC}$ is holomorphic with respect to the Hermitian
structure
$J$ on $U$ given by $J(q)=J(M(z(q))$ for all $q\in U$;

(ii)  this Hermitian structure is parallel along the fibres of $z(q)$,
i.e. the fibres of $z(q)$ are {\em superminimal\/} with respect to $J$;

(iii) $z(q)$ is a harmonic morphism;

(iv) $z(q)$ is submersive at points $q$ where $\frac{\partial F}{\partial
q} (q,z(q))\neq 0$.
\end{proposition}

\begin{remarks} (i)  This result also follows from \cite{Gu-Wo-1},
see\cite{Gu-Wo-2}.

 (ii) In the case $m=2$ this is a slight reformulation of the results  of
\cite{Wo-3}.  In that case it gives, locally, {\em all\/} submersive
harmonic morphisms
${\RR}^{4}\stackrel{open}{\supset}U\ra {\CC}$.  We shall see below that
 this is no longer the case when $m>2$.
\end{remarks}

Another special case of interest is when $k=m$ and $f$ is of the form
\begin{equation} f^{a}(w,z)=w^{a}-h^{a}(z) \label{eq:f-can}
\end{equation} with $h$ holomorphic. We shall call this the {\em
canonical form\/} of the equations (\ref{eq:f}).  $\big($The reason for
this name is that if the holomorphic function $f$ in Hypotheses
\ref{hyp:general} has Jacobian $(\pa f^a / \pa w^b)_{a,b = 1,\ldots m}$
non-singular, we can solve the equations (\ref{eq:f}) locally for $w =
h(z)$ and then the equations (\ref{eq:f}) can be written in the form
(\ref{eq:f-can}).$\big)$
  In this case we can calculate  the
$C^{ab}_{ij}$ in (\ref{eq:Lap2'}) by means of the following:

\begin{lemma} For any non-singular complex $k\times k$ matrix $K$, set
$A=K^{-1}$ and
$C^{ab}_{ij}=A^{a}_{i}A^{b}_{j}-A^{b}_{i}A^{a}_{j},\ (i,j,a,b=1, \ldots
,k)$. Then
\begin{equation} C^{ab}_{ij}=\frac{1}{\det K}(-1)^{i+j+a+b}\hbox{\rm
sign}(j-i)\hbox{\rm sign}(b-a)E_{ab}^{ij}
\label{eq:C2}
\end{equation} where for
$a<b,\ i<j$,\ $E_{ab}^{ij}$ is the determinant of the matrix obtained
from $K$ by omitting rows $i$ and $j$ and columns $a$ and $b$.
\end{lemma}

\proof  Since permuting the rows and columns of a determinant simply
multiplies it by the sign of the permutation, it suffices to prove the
lemma for $C^{12}_{12}$. Let $L=(L^{a}_{i})$ denote the $k\times
k$-matrix of cofactors  of $K$.   Then
$$ (\det K)^{2}C^{12}_{12}=L^{1}_{1}L^{2}_{2}-L^{1}_{2}L^{2}_{1}.
$$ Expanding the determinants $L^{a}_{i}$ along the top rows, we obtain
$$ (\det
K)^{2}C^{12}_{12}=\left\{\sum_{r=2}^{k}(-1)^{r}K^{r}_{2}E_{1r}^{12}\right\}
\left\{K^{1}_{1}E_{12}^{12}+\sum_{s=3}^{k}(-1)^{s}K^{s}_{1}E_{2s}^{12}\right\}
$$
$$-\left\{\sum_{r=2}^{k}(-1)^{r}K^{r}_{1}E_{1r}^{12}\right\}
\left\{K^{1}_{2}E_{12}^{12}+\sum_{s=3}^{k}(-1)^{s}K^{s}_{2}E_{2s}^{12}\right\}
$$
$$ =E_{12}^{12}\sum_{r=2}^{k}(-1)^{r}E_{1r}^{12}\left| \begin{array}{ll}
K^{1}_{1} & K^{r}_{1} \\ K^{1}_{2} & K^{r}_{2}
\end{array} \right|
+\sum_{r=2}^{k}\sum_{s=3}^{k}(-1)^{r+s}E_{1r}^{12}E_{2s}^{12}
\left| \begin{array}{ll} K^{s}_{1} & K^{r}_{1} \\ K^{s}_{2} & K^{r}_{2}
\end{array} \right|
$$
\begin{equation} =E_{12}^{12}\sum_{r=2}^{k}(-1)^{r}E_{1r}^{12}\left|
\begin{array}{ll} K^{1}_{1} & K^{r}_{1} \\ K^{1}_{2} & K^{r}_{2}
\end{array} \right|
-\sum_{\stackrel{r,s=2}{r<s}}^{k}(-1)^{r+s}(E_{1r}^{12}E_{2s}^{12} -
E_{1s}^{12}E_{2r}^{12})\left| \begin{array}{ll} K^{r}_{1} & K^{s}_{1} \\
K^{r}_{2} & K^{s}_{2}
\end{array} \right|. \label{eq:det}
\end{equation} The result will follow after we have established:
\eproof

\begin{claim}\label{th:Sylvester}
$$ E_{1r}^{12}E_{2s}^{12}-E_{1s}^{12}E_{2r}^{12}=E_{12}^{12}\hbox{\rm
sign}(s-r)E_{rs}^{12}.
$$
\end{claim}

For then the right hand side of (\ref{eq:det}) equals
$$ E_{12}^{12}\left\{\sum_{r=2}+{k}(-1)^{r}E_{1r}^{12}\left|
\begin{array}{ll} K^{1}_{1} & K^{r}_{1} \\ K^{1}_{2} & K^{r}_{2}
\end{array} \right|
-\sum_{\stackrel{r,s=2}{r<s}}^{k}(-1)^{r+s}E_{rs}^{12}\left|
\begin{array}{ll} K^{r}_{1} & K^{s}_{1} \\ K^{r}_{2} & K^{s}_{2}
\end{array} \right| \right\}
$$
$$ =-E_{12}^{12}\left\{\sum_{\stackrel{r,s=2}{r<s}}^{k}(-1)^{r+s}
E_{rs}^{12}\left| \begin{array}{ll} K^{r}_{1} & K^{s}_{1} \\ K^{r}_{2} &
K^{s}_{2}
\end{array} \right| \right\}.
$$ That the expression in brackets is $-\det K$ follows from Laplace's
Theorem  for determinants \cite{Metzler-Muir}.

\bigskip
\noindent {\bf Proof of Claim \ref{th:Sylvester}}\hskip 0.6em It suffices
to prove the case $r < s$.
 We apply Sylvester's Theorem \cite{Metzler-Muir}.  For this we need some
notation.  Denote the $(k-2)\times (k-2)$-determinant obtained from K  by
omitting  rows $1$ and $2$ and columns $r$ and $s$ by $\left| 12\cdots
(r-1)(r+1)
\cdots (s-1)(s+1)\cdots k\right|$.  Then
$$ E_{1r}^{12}E_{2s}^{12}\equiv \left| 23\cdots (r-1)(r+1)\cdots k\right|
\left| 13\cdots (s-1)(s+1)\cdots k\right|
$$ which by Sylvester's Theorem equals
$$
\left| 13\cdots (r-1)(r+1)\cdots k\right|\left| 23\cdots (s-1)(s+1)\cdots
k
\right|
$$
$$ +\left| 33\cdots (r-1)(r+1)\cdots k\right|\left| 124\cdots (s-1)(s+1)
\cdots k\right|
$$
$$ +\left| 43\cdots (r-1)(r+1)\cdots k\right|\left| 132\cdots (s-1)(s+1)
\cdots k\right|
$$
$$
\vdots
$$
$$   +\left| r3\cdots (r-1)(r+1)\cdots k\right|\left| 13\cdots (r-1)2(r+1)
\cdots (s-1)(s+1)\cdots k\right|
$$
$$
\vdots
$$
$$ =\left| 13\cdots (r-1)(r+1)\cdots k\right|\left| 23\cdots (s-1)(s+1)
\cdots k\right|
$$
$$ +\left| r3\cdots (r-1)(r+1)\cdots k\right|\left| 13\cdots (r-1)2(r+1)
\cdots (s-1)(s+1)\cdots k\right|,
$$ (all other terms vanishing due to their having two identical columns
in  a determinant).  But this is precisely
$E_{2r}^{12}E_{1s}^{12}+(-1)^{r+1} (-1)^{r+1}E_{12}^{12}E_{rs}^{12}$,
establishing the result.
\eproof

We can now find a formula for the Laplacian of a solution to
(\ref{eq:f-can}). Firstly the hypotheses:

\begin{hypotheses}\label{hyp:can} Suppose that
$h:{\Cm}\stackrel{open}{\supset}V\ra {\Cm}$ and $M:
{\Cm}\stackrel{open}{\supset}V\ra so(m,\CC)$ are given holomorphic
functions with $h$ of maximal rank.

Set
$$ F(q,z)=q-M(z)\bar{q}-h(z) \quad\hbox{for}\quad (q,z)\in {\Rmm}\times
V\ .
$$ We suppose that the Jacobian matrix
$K=(\partial F^{b}/\partial z^{a})_{a,b=1\cdots ,k}$ is non-singular on
some non-empty open set $A'\subset {\Rmm}\times V$.

Let $z : {\Rmm}
\stackrel{open}{\supset}U\ra V \subset {\Cm},\ z=z(q)$ be a smooth local
solution to equation (\ref{eq:F}):
\begin{equation} F(q,z) \equiv q - M(z)\bar{q} - h(z) = 0 \label{eq:F-can}
\end{equation}
\end{hypotheses}

Now the result:
\begin{proposition}\label{pr:rewrite} Suppose that Hypotheses
\ref{hyp:can} hold.   Then

(i)  $z:U\ra {\Ck}$ is holomorphic with respect to the Hermitian structure
$J$ on $U$ given by $J(q)=J(M(z(q))$ for $q\in U$;

(ii)  this Hermitian structure is parallel (i.e. constant) along the
fibres of
$z(q)$;

(iii)  the Laplacian of each component of $z(q)$ is given by
\begin{equation}
\frac{1}{4}\Delta z^a = \frac{(-1)^a}{\det K} \sum_{1 \leq i<j \leq m}
(-1)^{i+j-1}
\left| \begin{array}{ccccc}
\frac{\pa M^i_{\bar{j}}}{\pa z^1} & \ldots &
\widehat{\frac{\pa M^i_{\bar{j}}}{\pa z^a}} &
\ldots & \frac{\pa M^i_{\bar{j}}}{\pa z^m} \\ K^{k_1}_1 & \ldots &
K^{k_1}_a & \ldots & K^{k_1}_m \\ K^{k_{m-2}}_1 & \ldots & K^{k_{m-2}}_a
& \ldots & K^{k_{m-2}}_m \end{array}
\right| \label{eq:Lap-final}
\end{equation} where $(k_1,...k_{m-2}) = (1,\ldots , \hat{i} , \ldots
,\hat{j} , \ldots , m)$.
\end{proposition}

\proof  Parts (i) and (ii) follow from Proposition \ref{th:super}. For
part (iii), by Proposition \ref{th:super} we have that the Laplacian is
given by (\ref{eq:Lap2'}) where $A = K^{-1}$.  Then the $C^{ab}_{ij}$ are
given by (\ref{eq:C}).  To get  Formula (\ref{eq:Lap-final}) from
(\ref{eq:Lap2'}) and (\ref{eq:C}) simply note that $(-1)^b\hbox{\rm
sign}(b-a) = \hbox{\rm sign}\,\sigma$ where
$$
\sigma = \left( \begin{array}{ccccccccc} 1 & 2 & \ldots & \hat{a} &
\ldots & \ldots &
\ldots & \ldots & m \\ b & 1 & \ldots & \ldots & \hat{a} & \ldots &
\hat{b} & \dots & m
                                        \end{array} \right)
$$ and for $i < j$,
$$
\sum_{b=1, \dots , \hat{a}, \ldots , m}
                      \hbox{\rm sign}(\sigma) \frac{\pa M^i_{\bar j}}{\pa
z^b} E^{ij}_{ab}
$$ is the above determinant.
\eproof

We shall see that the data $h, M$ in Hypotheses \ref{hyp:can} can be
chosen such that  one or more of the components of  $z(q)$ is a harmonic
morphism.  In fact  this construction gives {\em all\/} harmonic
morphisms locally from  open subsets of $\Rmm$ to $\CC$ which are
holomorphic with respect to some Hermitian structure on $U$.  Precisely:

\begin{proposition}\label{prop:main} Let $\phi :\Rmm
\stackrel{open}{\supset} U \ra  {\CC}$ be holomophic with respect to some
Hermitian structure $J$ on
$U$.  (In particular, if $\phi$ is harmonic, it is a harmonic morphism.)

  Let $q_{0}\in U$ be a regular point of $\phi$.

(i)  By a permutation of the coordinates $(q^{i},\ q^{\bar{i}})$, we can
assume that $J(q_0)$ is  positive and lies in a large cell of
${\JJ}^{+}(\Rmm)$, so that $J(q_{0})=J(M)$  for some $M\in so(m,\CC)$.

(ii) Then in some neighbourhood $U' \subset U$ of $q_0$, the given map
$\phi$ is  the first component $z^{1}$ of a smooth solution $U'\ra\Cm $,\
$q\mapsto z(q)$
 to an equation
$$ F(q,z)=0
$$ with
$$ F(q,z)\equiv f(q-M(z)\bar{q}, z)
$$ where $f:{\Cm}\times {\Ck} \stackrel{open}{\supset} W \ra {\Ck}
$ and $M:{\Cm} \stackrel{open}{\subset} V \ra
so(m,\CC)={\CC}^{m(m-1)/2}$ are holomorphic maps and  $k$ is an integer,
$1\leq k \leq m$. The Laplacian of $\phi  = z^1$ is then given by
(\ref{eq:Lap2},
\ref{eq:A}) or, equivalently, (\ref{eq:Lap2'}, \ref{eq:C}).

(iii) We can choose $k=m$ and $f$ to be in canonical form
(\ref{eq:f-can}):
$$ f(w,z) = w - h(z)
$$ where $h: \Cm \stackrel{open}{\supset} V \ra {\Cm}$ is holomorphic and
of rank $m$  everywhere.  The Laplacian of $\phi$ is then given by
(\ref{eq:Lap-final}).
\end{proposition}

\proof   Since, on a neighbourhood of $q_{0}$, the map $\phi$ has maximal
rank (one), we can choose local complex coordinates $(z^{1},\cdots
,z^{m})$ on a possibly  smaller neighbourhood $U'$ of $q$ such that
$z^{1}=\phi$.  But then these  complex coordinates and the standard
complex coordinates $(w^{1},\cdots, w^{m})$ for $(U', J)$ are related by
a holomorphic diffeomorphism
$w=h(z)$ so that
$w-h(z)=0$ and we are done.
\eproof

\subsection{Low dimensional formulae} To find examples of harmonic
morphisms using Proposition (\ref{pr:rewrite}) it is useful write out
formula (\ref{eq:Lap-final}) in the cases $m=k=2 \hbox{ or } 3$.

\bigskip {\em The case $m=k=2$.}

Writing $\mu$ for $\mu^1$, the matrix $M$ reads
$$ M = M(\mu) =
\left(\begin{array}{rr} 0 & \mu\\ -\mu & 0
\end{array}\right)
$$ so that the standard complex coordinates for the Hermitian structure
$J = J(M)$ (see
\S 2) are given by
$$
\left.\begin{array}{ll} w^{1}=q^{1}-\mu q^{\bar{2}}\\ w^{2}=q^{2}+\mu
q^{\bar{1}}
\end{array}\right\}.
$$   Equations (\ref{eq:F}) for $z:q = (q^1,q^2) \mapsto z = (z^1,z^2)$
read
$$ f(q^{1}-\mu(z)q^{\bar{2}}, q^{2}+\mu(z)q^{\bar{1}}, z) = 0.
\label{eq:f-4d}
$$  When $f$ is in canonical form (\ref{eq:f-can}) these read
\begin{equation}\label{eq:can2d}
\left.\begin{array}{ll} F^1(q,z) \equiv q^1 - \mu(z^1,z^2)q^{\bar{2}} -
h^1(z^1,z^2) =  0, \\ F^2(q,z) \equiv q^2 + \mu(z^1,z^2)q^{\bar{1}} -
h^2(z^1,z^2) =  0.
\end{array}\right\}
\end{equation} Then (\ref{eq:Lap-final}) simplifies to give:

\begin{proposition} Suppose Hypotheses \ref{hyp:can} hold with $m=2$.
Then the Laplacians of the components of $z(q)$ are given by:

\begin{eqnarray*}
\frac{1}{4}\Delta z^{1} & = & \hskip 0.7em \frac{1}{\det K}\frac
{\partial \mu}{\partial z^{2}}\ , \\
\frac{1}{4}\Delta z^{2} & = & -\frac{1}{\det K}\frac{\partial\mu}
{\partial z^{1}}.
\end{eqnarray*}
\end{proposition} Since by Theorem \ref{prop:main}, any submersive
holomorphic map is locally the solution to equations (\ref{eq:can2d}),
these formulae confirm the result of
\cite{Wo-3} that a submersive harmonic  map ${\RR}^{4}\supset U\ra {\CC}$
which is holomorphic with respect to a Hermitian  structure on $U$ is a
harmonic morphism if and only if its fibres are  superminimal.

\bigskip {\em The case $m=k=3$.}

 This time the matrix $M$ reads
\begin{equation} M = M(\mu_1,\mu_2,\mu_3) =
\left(\begin{array}{rrr} 0 & \mu _{1} & \mu _{2}\\ -\mu _{1} & 0 & \mu
_{3}\\ -\mu _{2} & -\mu _{3} & 0
\end{array}\right),
\label{eq:M3} \end{equation} so that the standard complex coordinates for
the Hermitian structure $J = J(M)$ are given by
$$
\left.\begin{array}{lll} w^{1}=q^{1}-\mu _{1}q^{\bar{2}}-\mu
_{2}q^{\bar{3}}\\ w^{2}=q^{2}+\mu _{1}q^{\bar{1}}-\mu _{3}q^{\bar{3}}\\
w^{3}=q^{3}+\mu _{2}q^{\bar{1}}+\mu _{3}q^{\bar{2}}
\end{array}\right\}.
$$ With $f$ in canonical form (\ref{eq:f-can}) the equations
(\ref{eq:F-can}) for $z :q = (q^1,q^2,q^3) \mapsto z = (z^1,z^2,z^3)$ read
\begin{equation}
\left.\begin{array}{lll} F^{1}(q,z)\equiv q^{1}-\mu
_{1}(z)q^{\bar{2}}-\mu _{2}(z)q^{\bar{3}} -h^{1}(z)=0\\ F^{2}(q,z)\equiv
q^{2}+\mu _{1}(z)q^{\bar{1}}-\mu _{3}(z)q^{\bar{3}} -h^{2}(z)=0\\
F^{3}(q,z)\equiv q^{3}+\mu _{2}(z)q^{\bar{1}}+\mu _{3}(z)q^{\bar{2}}
-h^{3}(z)=0
\end{array}\right\}\label{eq:F3}
\end{equation} and  from equations (\ref{eq:Lap-final}) we obtain

\begin{proposition} Suppose Hypotheses \ref{hyp:can} hold with $m=3$.
Then the Laplacians of the components of $z(q)$ are given by:
\begin{eqnarray}
\frac{1}{4}\Delta z^{1} & = & \frac{1}{\det
K}\left\{\left|\begin{array}{ll}
\frac{\partial \mu _{1}}{\partial z^{2}} &
\frac{\partial \mu _{1}}{\partial z^{3}}  \\ K^{3}_{2} &
K^{3}_{3}\end{array}\right|  -\left|\begin{array}{ll}\frac{\partial \mu
_{2}}{\partial z^{2}} &
\frac{\partial \mu _{2}}{\partial z^{3}}\\ K^{2}_{2} &
K^{2}_{3}\end{array}\right| +\left|\begin{array}{ll}\frac{\partial \mu
_{3}}{\partial z^{2}} &
\frac{\partial \mu _{3}}{\partial z^{3}} \\ K^{1}_{2} &
K^{1}_{3}\end{array}\right|\right\} \nonumber \\
  \frac{1}{4}\Delta z^{2} & = & \frac{1}{\det
K}\left\{\left|\begin{array}{ll}
\frac{\partial \mu _{1}}{\partial z^{3}} &
\frac{\partial \mu _{1}}{\partial z^{1}}\\ K^{3}_{3} &
K^{3}_{1}\end{array}\right|  -\left|\begin{array}{ll}\frac{\partial \mu
_{2}}{\partial z^{3}} &
\frac{\partial \mu _{2}}{\partial z^{1}}\\ K^{2}_{3} &
K^{2}_{1}\end{array}\right| +\left|\begin{array}{ll}\frac{\partial \mu
_{3}}{\partial z^{3}} &
\frac{\partial \mu _{3}}{\partial z^{1}}\\ K^{1}_{3} &
K^{1}_{1}\end{array}\right|\right\} \label{eq:Lap3ex} \\
  \frac{1}{4}\Delta z^{3} & = & \frac{1}{\det
K}\left\{\left|\begin{array}{ll}
\frac{\partial \mu _{1}}{\partial z^{1}} &
\frac{\partial \mu _{1}}{\partial z^{2}}\\ K^{3}_{1} &
K^{3}_{2}\end{array}\right|  -\left|\begin{array}{ll}\frac{\partial \mu
_{2}}{\partial z^{1}} &
\frac{\partial \mu _{2}}{\partial z^{2}}\\ K^{2}_{1} &
K^{2}_{2}\end{array}\right| +\left|\begin{array}{ll}\frac{\partial \mu
_{3}}{\partial z^{1}} &
\frac{\partial \mu _{3}}{\partial z^{2}} \\ K^{1}_{1} &
K^{1}_{2}\end{array}\right|\right\}\nonumber
\end{eqnarray} where
\begin{eqnarray} K^{1}_{a}=\frac{\partial F^{1}}{\partial z^{a}}=
-\frac{\partial \mu _{1}}{\partial z^{a}}q^{\bar{2}} -\frac{\partial \mu
_{2}}{\partial z^{a}}q^{\bar{3}} -\frac{\partial h^{1}}{\partial z^{a}}
\nonumber \\ K^{2}_{a}= \frac{\partial F^{2}}{\partial z^{a}} = \hskip
0.7em
\frac{\partial \mu _{1}}{\partial z^{a}}q^{\bar{1}} -\frac{\partial \mu
_{3}}{\partial z^{a}}q^{\bar{3}} -\frac{\partial h^{2}}{\partial z^{a}}\\
K^{3}_{a}= \frac{\partial F^{3}}{\partial z^{a}} = \hskip 0.7em
\frac{\partial \mu _{2}}{\partial z^{a}}q^{\bar{1}} +\frac{\partial \mu
_{3}}{\partial z^{a}}q^{\bar{2}} -\frac{\partial h^{3}}{\partial z^{a}}
\nonumber
\label{eq:K3}
\end{eqnarray}
\end{proposition}

\section{Examples and results}\label{sec:examples}

In this section we shall construct some examples using the theory of the
last section together with Proposition \ref{prop:sethaho}. The idea is to
choose the data $h$ and $M=M(\mu)$ in Hypotheses \ref{hyp:can} and to
solve Equation (\ref{eq:F-can}) for $z = \phi(q)$.  Then the Laplacians
of the components of $z$ are found by Equation (\ref{eq:Lap-final}) and
the data is chosen judiciously such that one or more of the components is
a harmonic morphism with specific properties.  Sometimes, $z$ is composed
with a futher holomorphic map $\phi$ to obtain more harmonic holomorphic
maps by Proposition \ref{prop:sethaho}. Before doing this, we must
consider, in the first subsection, three important properties of such
examples.  Readers wishing to get straight to the examples should skip to
\S
\ref{subsec:examples}.

\subsection{Fullness, superminimality and K\"ahler structures}
\label{subsec:fullsuperk}

{\em 1. Fullness and reduction}

\begin{definition}\label{def:full} Call a map $\phi:\RR^n
\stackrel{open}{\supset} U
\ra \CC$ {\em full\/} if we cannot write it as $\phi = \psi \circ \pi_A$
for some orthogonal projection $\pi_A$ onto a subspace $A$ of $\RR^n$ and
map $\psi:\pi_A(U)
\ra
\CC$.  If, on the other hand $\phi$ does so factor we say that $\phi$
{\em reduces\/} to $A$ and $\psi$ is a {\em reduction\/} of
$\phi$.
\end{definition}

 Note that if $\phi$ does so factor, then it is a harmonic map (resp.
harmonic morphism) if and only if $\psi$ is a harmonic map (resp.
harmonic morphism), see
\cite{Gu-3}.

Similarly, we shall say that a foliation ${\cal F}$ on $U$ is {\em
full\/} if it is not the inverse image $(\pi_A)^{-1}({\cal F}^{\prime})$
of a foliation ${\cal F}^{\prime}$ on $\pi_A(U)$ for any subspace $A$.
\smallskip

We now give a test for fullness for holomorphic maps $\phi\circ z$ where
$z(q)$ is given by a solution of equation (\ref{eq:F-can}). To be precise
we assume Hypotheses
\ref{hyp:can} and that $\phi:V \ra \CC$ is a given holomorphic map. We
seek conditions that $\phi \circ z$ factor to $\RR^{2m-1}$.  As usual,
write
$K^a_b = \pa F^a/\pa z^b$,\ $(a,b=1,\ldots m)$ and $(\hat{K}^a_b)$ for
the matrix of cofactors of $K$ so that $(K^{-1})^a_b = (1/\det K)
\hat{K}^b_a$.

Let $A$ be a $(2m-1)$-dimensional subspace and let $v = a^j \pa/\pa x^j$
be a non-zero vector orthogonal to it.  Define
$\alpha \in \Cm$ by $\alpha^j = a^{2j-1}-\ii a^{2j}$,\
$(j=1,\ldots m)$.
\begin{proposition}\label{prop:fulltest}
 Suppose Hypotheses \ref{hyp:can}, and let $\phi:V \ra \CC$ be a
holomorphic function.  Then the map
$U \ra \CC$,\ $q \mapsto \phi(z(q))$ factors to $\pi_A(U)\subset A$ if
and only if
\begin{equation}
\sum_{a=1}^m \frac{\pa\phi}{\pa z^a} \left| \begin{array}{lllllll} K^1_1 &
\cdots & K^1_{a-1} & w^1(\alpha,z(q)) & K^1_{a+1} & \cdots & K^1_m \\
K^2_1 &
\cdots & K^2_{a-1} & w^2(\alpha,z(q)) & K^2_{a+1} & \cdots & K^2_m \\
\cdots & \cdots & \cdots   & \cdots      & \cdots    & \cdots & \cdots \\
K^m_1 & \cdots & K^m_{a-1} & w^m(\alpha,z(q)) & K^m_{a+1} & \cdots & K^m_m
\end{array}
\right| = 0 \hbox{ for all } q\in U \label{eq:redn-Eucl}
\end{equation} where the $K^a_b$ are evaluated at $(q,z(q))$.
\end{proposition}

\proof A function $\psi:U \ra \CC$ factors if and only if its directional
derivative in direction $v$ is zero, i.e.
\begin{equation}
\alpha^I \pa \psi/\pa q^I = 0 \label{eq:deriv-Eucl}
\end{equation} (where we sum over $I = 1,\ldots ,m,\bar{1},\ldots
\bar{m}$). Now, setting $\psi = \phi \circ z$, (\ref{eq:dzdq}) gives
\begin{equation}
\frac{\pa \psi}{\pa q^I} = \frac{\pa \phi}{\pa z^a} \frac{\pa z^a}{\pa
q^I} = - \frac{\pa \phi}{\pa z^a} (K^{-1})^a_b \frac{\pa F^b}{\pa q^I} =
- \frac{\pa
\phi}{\pa z^a} \frac{1}{\det K} \hat{K}^b_a \frac{\pa F^b}{\pa q^I}
\label{eq:dzdq2}
\end{equation} so that (\ref{eq:deriv-Eucl}) reads
$$
\frac{\pa\phi}{\pa z^a}(\alpha^I \frac{\pa F^b}{\pa q^I})
\hat{K}^b_a = 0.
$$ Noting that $\alpha^I \pa F^b/\pa q^I =
\alpha^I\pa w^b/\pa q^I = w^b(\alpha, z)$ (the last equality using the
fact that
$w(\alpha ,z) = \alpha - M\bar{\alpha}$ is homogeneous of degree one in
$\alpha$), we obtain the formula (\ref{eq:redn-Eucl}).
\eproof

\begin{example}  Putting $a = (1,0,\ldots,0)$ we see that the function
$z^1$ factors through the natural projection $\Rmm\ra\RR^{2m-1}$,\
$(x^1,x^2, \ldots ,x^{2m})\mapsto (x^2,\ldots ,x^{2m})$, i.e. is
independent of $x^1$, if and only if

\begin{equation}
\left|\begin{array}{llll} 1       & K^1_2 & \ldots & K^1_m \\
\mu_{1} & K^2_2 & \ldots & K^2_m \\
\ldots  & \ldots & \ldots & \ldots \\
\mu_{m-1} & K^m_2 & \ldots & K^m_m
\end{array}\right| = 0 \quad\hbox{for all}\quad q\in U
\label{eq:redn-Eucl0}
\end{equation} where the $K^a_b$ are evaluated at $(q,z(q))$.
\end{example}

\bigskip {\em 2. Superminimality}

The definition of superminimality was given in Definition
\ref{def:superminimal}. It is easy to show that if an even dimensional
submanifold
$S$ of an open subset of $\Rmm$ equipped with an adapted almost Hermitian
structure
$J$ is superminimal with respect to $J$ then its Weingarten map $A$
satisfies
\begin{equation} A_{JX}V = JA_X V \quad \hbox{for all } x \in S,\  X \in
T_x S \hbox{ and }  V \in T_x\Rmm \hbox { normal to } M\ ;
\label{eq:super-Wein}
\end{equation}  from this it easily follows that an even dimensional
submanifold which is superminimal with respect to an almost Hermitian
structure is minimal.  We further note that if $S$ is of real codimension
$2$, the property (\ref{eq:super-Wein}) characterizes superminimality
with respect to an almost Hermitian structure (see, for example,
\cite{Gu-Wo-2}).

We can characterize this condition geometrically by
\begin{proposition} A codimension $2$ submanifold $S$ of $\Rmm$ is
superminimal with respect to some adapted almost Hermitian structure if
and only if, for all
$x
\in S$,\
$A_Z V$ is isotropic for all $V \in T_x S$ and all isotropic $Z \in
T_x\Rmm$ normal to
$S$.
\end{proposition}

We also give a criterion for supermininality with respect to an almost
complex structure:
\begin{proposition}\label{prop:supcrit} Let $\psi:\Rmm
\stackrel{open}{\supset} U \ra \CC$ be a submersive harmonic morphism (or
any submersion which satisfies Equation (\ref{eq:HWC})).  Then a
connected component of a fibre $E$ of $\psi$ is superminimal with respect
to some adapted almost complex structure $J$ on $U$ if

(1) there exists an $m$-dimensional complex subspace $W'$ of $\Cmm$ with
$W' \cap \bar{W'} = \{0\}$ and
\begin{equation}
\nabla\psi = (\pa\psi/\pa x^1, \ldots ,\pa\psi/\pa x^{2m}) \in W'
\label{eq:gradient}
\end{equation} at all points $q$ of $E$, or equivalently,

(2) there exists an $m$-dimensional complex subspace $W$ of $\Cmm$ with
$W \cap
\bar{W} = \{0\}$ such that
\begin{equation}
 b^I \frac{\pa\psi}{\pa q^I} = 0 \label{eq:superfull}
\end{equation} for all $b \in W,\ q \in E$.
\end{proposition}

\proof For the ``only if" part of (1), set $W' =$ the $(1,0)$-tangent
space of $\Rmm$ with respect to $J$.  This satisfies the desired
conditions.   For the converse, choose $J$ to have $+\ii$- (resp.
$-\ii$-) eigenspace $W'$ (resp. $W$) at all points of $E$. To see the
equivalence of (1) and (2), set $W =$ the orthogonal complement of
$W'$ in $\Cmm$ noting that (\ref{eq:superfull}) says that $\nabla\psi$ is
orthogonal to $W$.
\eproof

For the fibres of a map $\phi\circ z$  where $z(q)$ is a solution to
(\ref{eq:F-can}) we deduce the following test for superminimality:

\begin{proposition}\label{prop:suptest} Suppose Hypotheses
\ref{hyp:can}.  Additionally, let $\phi:V \ra \CC$ be a holomorphic
function. Then a connected component $E$ of a fibre of $q \mapsto
\phi(z(q))$ is superminimal with respect to some almost complex structure
if and only if there is an $m$-dimensional linear subspace $W$ of $\Cmm$
with $W \cap \bar{W} = 0$ such that
\begin{equation}
\sum_{a=1}^m \frac{\pa\phi}{\pa z^a} \left| \begin{array}{lllllll} K^1_1
& \cdots & K^1_{a-1} & w^1(b,z(q)) & K^1_{a+1} & \cdots & K^1_m \\ K^2_1
& \cdots & K^2_{a-1} & w^2(b,z(q)) & K^2_{a+1} & \cdots & K^2_m \\
\cdots & \cdots & \cdots   & \cdots      & \cdots    & \cdots & \cdots \\
K^m_1 & \cdots & K^m_{a-1} & w^m(b,z(q)) & K^m_{a+1} & \cdots & K^m_m
\end{array}
\right| = 0 \label{eq:suptest}
\end{equation} holds for all $b \in W$ and all $q \in E$.
\end{proposition}

\proof Apply the last Proposition to $\psi = \phi \circ z$ calculating
its gradient as in Proposition \ref{prop:fulltest}.
\eproof

\bigskip {\em 3. K\"ahler structures}

Lastly we consider whether a given map $\psi = \phi \circ z:\Rmm
\stackrel{open}{\supset} U \ra \CC$ is holomorphic with respect to a
parallel complex structure on $\Rmm$, for example a K\"ahler structure.
This is very similar to the last test, only now we insist that the linear
subspace does not vary from fibre to fibre.

\begin{proposition}\label{prop:kahlertest} Suppose Hypotheses
\ref{hyp:can}. Additionally, let $\phi:V \ra \CC$ be a holomorphic
function. Then
$q \mapsto \phi(z(q))$ is holomorphic with respect to a parallel complex
structure on
$\Rmm$ if and only if there is an $m$-dimensional linear subspace $W$ of
$\Cmm$ with $W \cap W' = {0}$ such that Equation (\ref{eq:suptest}) holds
for all
$b \in W$ and $q \in U$.
\end{proposition}

To include the odd dimensional case we introduce some terminology:
\begin{definition}\label{def:arise} Say that a map
$\phi:\RR^n \stackrel{open}{\supset} U \ra \CC$ {\em arises from a
K\"ahler structure\/} if either
$n$ is even and it is holomorphic with respect to a K\"ahler structure on
$\RR^n$, or
$n$ is odd and it is the reduction of a map on an open subset of
$\RR^{n+1}$ which is holomorphic with respect to a K\"ahler structure on
$\RR^{n+1}$.
\end{definition}

Note that in the second case ($n$ odd), $\phi$ further reduces to a map
on an even dimensional subspace of $\RR^n$ (for it is invariant in the
direction $v \in {(\RR^n)}^{\perp} \cap \RR^{n+1}$ and so it is also
invariant in the direction $Jv$.  Hence it reduces to the orthogonal
complement of $Jv$ in
$\RR^n$.

\subsection{Examples of harmonic morphisms on even dimensional Euclidean
spaces}
\label{subsec:examples}

In this subsection we construct examples of harmonic morphisms from open
subsets of
${\RR}^{2m}$ by using Propositions \ref{pr:rewrite} and
\ref{prop:sethaho}.

\bigskip {\em Examples with superminimal fibres with respect to a
Hermitian structure which is not cosymplectic}

Suppose that in Proposition \ref{pr:rewrite} we choose
$\mu _{i}:{\Cm}
\supset V\ra {\CC}$ holomorphic with $\mu _{i}=\mu _{i}(z^{1}), i=1,
\cdots ,m(m-1)/2$, depending on $z^{1}$ only, with the holomorphic
functions $h^{i}$  chosen arbitrarily.  Then by (\ref{eq:Lap-final}),
$\Delta z^{1}=0$.  In this case  the Hermitian structure $J(M(\mu))$  is
constant along the fibres of $z^{1}$ which are therefore superminimal
with respect to it.

\begin{example}\label{ex:LSnotC} Let $\mu _1=z^1, \mu _2=z^1, \mu _3=0,
h^1=z^2, h^2=z^3, h^3=z^2z^3$ .  Then $z^1$ satisfies the system of
equations
$$
\left. \begin{array}{llll} q^1 - z^1q^{\bar{2}} & -z^1q^{\bar{3}} & -z^2
& =0 \\ q^2 + z^1q^{\bar{1}} &                & -z^3 & =0 \\ q^3 +
z^1q^{\bar{1}} &                & -z^2z^3 & =0
\end{array}
\right\}\ .
$$ Then $\det K$ is non-zero almost everywhere and, eliminating $z^2$ and
$z^3$ we see that
$z^1$ is determined by the quadratic equation
$$ (z^1)^2q^{\bar{1}}(q^{\bar{2}}+q^{\bar{3}}) + z^1\left[ q^{\bar{1}} +
q^2(q^{\bar{2}}+q^{\bar{3}}) - |q^1|^2\right] + q^3-q^1 q^2 = 0.
$$ Neither $z^2$ nor $z^3$ is harmonic hence the corresponding complex
structure is not cosymplectic.  Furthermore, it can be easily checked
that $z^1$ is not holomorphic with respect to any K\"ahler structure.

This example generalises to all $m>2$ by setting $h^1=z^2, h^2=z^3,
\ldots , h^{m-1}=z^m, h^m=z^2z^3\cdots z^m$, with the $\mu _i$'s
appropriately chosen functions of $z^1$.
\end{example}

{\em Examples with superminimal fibres with respect to a Hermitian
structure which is cosymplectic}

Note firstly that
 for $m=2$,\ $J$ cosymplectic implies $J$ K\"ahler; this is no longer the
case  for $m > 2$.  Suppose again that  $\mu _i=\mu _i(z^1)$ for each
$i=1,2,\ldots ,m(m-1)/2$, so that
$z^1$ has superminimal fibres with  respect to $J=J(M(\mu ))$, so that. In
particular $z^1$ is harmonic and by Proposition
\ref{prop:cosymplectic}, $J$ is cosymplectic if and only if
$z^2,z^3,\ldots z^m$ are also harmonic.  For simplicity of exposition we
consider  the case $m=3$.

By equation (\ref{eq:Lap3ex})
\begin{eqnarray*}
\frac{\det K}{4}\Delta z^2 & = & \frac{\partial\mu _1}{\partial
z^1}\frac{\partial h^3} {\partial z^3}-\frac{\partial\mu _2}{\partial
z^1}\frac{\partial h^2} {\partial z^3}  +\frac{\partial\mu _3}{\partial
z^1}\frac{\partial h^1} {\partial z^3} \\
  & = & \frac{\partial}{\partial z^3}(h^3\frac{\partial \mu _1}{\partial
z^1} -h^2\frac{\partial \mu _2}{\partial z^1}+h^1\frac{\partial \mu
_3}{\partial z^1})
\end{eqnarray*} Similarly

\begin{equation}
\frac{\det K}{4}\Delta z^3=
\frac{\partial}{\partial z^2}(h^3\frac{\partial \mu _1}{\partial z^1}
-h^2\frac{\partial \mu _2}{\partial z^1}+h^1\frac{\partial \mu
_3}{\partial z^1})
\end{equation} so that $J$ is cosymplectic if and only if

\begin{equation} h^3\frac{\partial \mu _1}{\partial z^1}
-h^2\frac{\partial \mu _2}{\partial z^1}+h^1\frac{\partial \mu
_3}{\partial z^1}=\alpha (z^1),
\end{equation} where $\alpha$ is a holomorphic function of $z^1$ only.
Then $z^1$ satisfies the single equation

\begin{equation}
\frac{\partial\mu _1}{\partial z^1}(q^3+\mu _2q^{\bar{1}}+\mu
_3q^{\bar{2}})  -\frac{\partial\mu _2}{\partial z^1}(q^2+\mu
_1q^{\bar{1}}-\mu _3q^{\bar{3}}) +\frac{\partial\mu _3}{\partial
z^1}(q^1-\mu _1q^{\bar{2}}-\mu _2q^{\bar{3}}) =
\alpha(z^1 ).
\end{equation}

Conversely, given $\mu _1,\mu _2,\mu _3$ as functions of $z^1$ and a
solution $z^1$ to the above equation,  setting $h^1=q^1-\mu
_1q^{\bar{2}}-\mu _2q^{\bar{3}}, h^2=q^2+\mu _1q^{\bar{1}}-\mu
_3q^{\bar{3}}, h^3=q^3+\mu _2q^{\bar{1}}+\mu _3q^{\bar{2}}$, we see that
$z^1$ is holomorphic with respect to a cosymplectic structure.

We can now easily generate non-trivial solutions.  For instance

\begin{example}\label{ex:LSC} Set $\mu _1=z^1, \mu _2=(z^1)^2, \mu
_3=(z^1)^3, h^3=2z^1z^2 - 3(z^1)^2 z^3, h^2=z^2, h^1=z^3$.  Then the
above procedure gives a harmonic morphism $z^1$ of the type considered
given by any local solution to the equation
$$ (z^1)^4q^{\bar{3}} + 2(z^1)^3q^{\bar{2}} + z^2(q^{\bar{1}} - 3q^1) +
2z^1q^2 - q^3 = 0.
$$

\bigskip {\em Examples with non-superminimal fibres}

Choose $\mu _i:{\Cm}\supset V\ra {\CC}$ holomorphic with $\mu _1=\mu
_1(z^1,\ldots , z^{m-1}), \mu _a=\mu _a(z^1),\ (a=2,\ldots , m(m-1)/2)$.
Then
\begin{equation}
\Delta z^1=-\frac{1}{\det K}
\left| \begin{array}{lllll}
\frac{\partial \mu _1}{\partial z^2} & \frac{\partial \mu _1}{\partial
z^3} & \cdots & \frac{\partial \mu _1}{\partial z^{m-1}} & 0 \\ K^3_2 &
K^3_3 & \cdots & K^3_{m-1} & K^3_m \\ K^4_2 & K^4_3 & \cdots & K^4_{m-1}
& K^4_m \\
\vdots & \vdots & \vdots & \vdots & \vdots \\ K^m_2 & K^m_3 & \cdots &
K^m_{m-1} & K^m_m
\end{array}\right| ,
\end{equation} with $K^a_b=-\partial h^a/\partial z^b,\ (a,b\geq 2$).

Then if $h^3, h^4, \ldots ,h^m$ are independent of $z^4$ the last column
of the above determinant vanishes and so $\Delta z^1=0$.  Noting that
unless $\mu _1$ depends on $z^1$ only, the complex structure $J=J(\mu )$
varies with $z^2, \cdots , z^{m-1}$ as well as with $z^1$, giving local
examples with fibres not superminimal with respect to $J$.  However it
might still be the case that
$z^1$ is holomorphic with respect to another Hermitian structure with
respect to which it has superminimal fibres, as the following example
shows.
\end{example}

\begin{example}\label{ex:LNSNC} Let $\mu _{1}=z^{2}, \mu _{2}=z^{1}, \mu
_{3}=0,  h^{1}=z^{3}, h^{2}=z^{3}, h^{3}=z^{2}$.  Then
$$
\det K=(q^{\bar{1}})^{2}+q^{\bar{1}}q^{\bar{2}}+q^{\bar{3}}
$$ is non-zero on the complement $U$ of the surface
$(q^{\bar{1}})^{2}+q^{\bar{1}}q^{\bar{2}}+q^{\bar{3}}=0$ and
$$ z^{1}=\frac{-q^{3}(q^{\bar{1}}+q^{\bar{2}})+q^{1}-q^{2}}
{(q^{\bar{1}})^{2}+q^{\bar{1}}q^{\bar{2}}+q^{\bar{3}}}
$$  is a harmonic morphism from $U$ to $\CC$, holomorphic wih respect  to
the Hermitian structure represented by
$$ M= M(z(q)) = \left(\begin{array}{rrr} 0 & z^{2} & z^{1}\\ -z^{2} & 0 &
0\\ -z^{1} & 0 & 0
\end{array}\right)
$$

 Furthermore $\Delta z^{2}=-4/\det K$  and $\Delta z^3=
-2q^{\bar{1}}/\det K$ which are non-zero, so that $J$ is not
cosymplectic.  However a simple check, using Theorem \ref{prop:supcrit}
shows that along each connected fibre component,
$\nabla z^1$ is contained in a fixed ${\CC}^3\subset {\CC}^6$, showing
that in fact the fibres of $z^1$ are superminimal with respect to some
complex  structure
\end{example}

This illustrates the possibility of a harmonic morphism being holomorphic
with respect to more than one Hermitian structure.  In fact there are
examples which are holomorphic with respect to a family of Hermitian
structures:

\begin{example}\label{ex:family} Let $m=3$ and suppose $\mu _1=\mu
_1(z^1,z^2), h^3=h^3(z^1)$.  Then the harmonic morphism defined by
equations (\ref{eq:F3}) is determined by the last of these equations:
$$ q^3+\mu _2(z^1)q^{\bar{1}}+\mu _3(z^1)q^{\bar{2}}-h^3(z^1)=0
$$ and provided $\mu _2, \mu _3$ and $h^3$ are chosen judiciously, for
example,
$h^3=z^1, \mu _2=(z^1)^2, \mu _3=(z^1)^3$, the resulting map is full.

We are now at liberty to choose $\mu _1, h^1, h^2$ as we wish, provided
of course that the system gives well-defined solutions.  To be specific,
if we set
$h^1=z^2, h^2=z^3, \mu _1=\mu _1(z^1, z^2)$, then the Jacobian matrix $K$
has determinant
$$
\left(\frac{\partial \mu _2}{\partial z^2}q^{\bar{2}}-1\right)
(2z^1q^{\bar{1}}+3(z^1)^2q^{\bar{2}}-1),
$$ which is non-zero in general, giving a family of Hermitian structures
param\-et\-ris\-ed by the function
$\mu _1(z^1, z^2)$ with respect to each of which $z^1$ is holomorphic.
Amongst this family are those parametrised by functions $\mu _1$
independent of $z^2$ and, with respect to such a Hermitian structure,
$z^1$ has superminimal fibres.  There are examples of harmonic
holomorphic maps which do not have superminimal fibres (with respect to
any almost complex structure) which we shall describe shortly when we
consider global examples.

The above construction generalizes to higher dimensions.  For example, if
$m=4$, set
$\mu _2=z^1, \mu _3=\mu _4=0, \mu _5=z^1, \mu _6=(z^1)^3, h^3=z^1z^2,
h^4=z^2$, and suppose $\mu _1=\mu _1(z^1, z^2)$.  Then $z^1$ is harmonic
and is determined by the last two equations of the system
(\ref{eq:F-can}):
$$
\left. \begin{array}{lll}

q^1-\mu _1q^{\bar{2}}-\mu _2q^{\bar{3}}-\mu _3q^{\bar{4}}-h^1 & = & 0 \\
q^2+\mu _1q^{\bar{1}}-\mu _4q^{\bar{3}}-\mu _5q^{\bar{4}}-h^2 & = & 0 \\
q^3+\mu _2q^{\bar{1}}+\mu _4q^{\bar{2}}-\mu _6q^{\bar{4}}-h^3 & = & 0 \\
q^4+\mu _3q^{\bar{1}}+\mu _5q^{\bar{2}}+\mu _6q^{\bar{3}}-h^4 & = & 0 .
\end{array} \right\}
$$

Eliminating $z^2$ from these gives $z^1$ as a solution of
$$
(z^1)^4q^{\bar{3}}+(z^1)^3q^{\bar{4}}+(z^1)^2q^{\bar{2}}+
z^1(q^4-q^{\bar{1}})-q^3=0.
$$ Then $z^1$ is full and is holomorphic with respect to Hermitian
structures param\-et\-ris\-ed by the function$\mu _1=\mu _1(z^1, z^2)$.
This example is additionally interesting because $z^1$ is invariant under
radial translation
$q\mapsto aq, 0<a<\infty $, and so reduces to a full harmonic morphism on
a domain in $S^7$ (see
\cite{Ba-Wo-6}).  Generalisation of these constructions to arbitrary
$m>2$ is straightforward.
\end{example}

{\em Global examples in even dimensions}

The above examples show the richness and variety that abounds for $m>2$.
We may ask whether such diversity can occur if we insist our examples be
globally defined on
${\Rmm}$.  We now turn our attention to this case.

Let $m>2$ and suppose $\mu _1=\mu _1(z^3,\ldots ,z^m),
\mu _2=\mu _2(z^4,\ldots ,z^m), \ldots ,\mu _{m-2}=\mu _{m-2}(z^m)$ with
$\mu _{m-1}=a_{m-1},\ldots ,\mu _{m(m-1)/2}=a_{m(m-1)/2}$ arbitrary
constants.  Define $h^k$ to be of the form $h^k=z^k+g_k(z^{k+1}, \ldots
,z^m)$, for each
$k=1,\ldots ,m-1$ and $h^m = z^m$, where the $g_k$ are arbitrary globally
defined holomorphic functions.  Then $M$ has the form
\begin{equation} M=\left(\begin{array}{rrrrrr} 0 & \mu _1 & \mu _2 &
\cdots & \mu _{m-2} & a_{m-1} \\ -\mu _1 & 0 & a_m & \cdots & \vdots &
\vdots \\ -\mu _2 & -a_m & \vdots & \cdots & \vdots & \vdots \\
\vdots & \vdots & \vdots & \cdots & \vdots & \vdots \\ -\mu _{m-2} &
\vdots & \vdots & \cdots & \vdots & \vdots \\ -a_{m-1} & \cdots & \cdots
& \cdots & \cdots &
\end{array}\right)
\end{equation} giving a family of cosymplectic structures defined
globally on ${\Rmm}$.  In fact
\begin{equation} K=\left(\begin{array}{rrrrr} -1 & * & * & \cdots & * \\
0 & -1 & * & \cdots & * \\ 0 & 0 & -1 & \cdots & * \\
\vdots & \vdots & \vdots & \cdots & * \\  0 & 0 & 0 & \cdots & -1
\end{array}\right)
\end{equation} has determinant $(-1)^m$ everywhere.  The system
(\ref{eq:F-can}) determines globally defined harmonic holomorphic
functions $z^1, \ldots , z^m$ which are found explicitly by first solving
the $m$th equation for $z^m$, then substituting into the $(m-1)$st to
obtain $z^{m-1}$ etc.  We fix our attention on a special case.

\begin{example}\label{ex:glob1}

Let $\mu _1=z^m, \mu _2=\cdots =\mu _{m(m-1)/2}=0, h^1=z^1, h^2=z^2,
\ldots , h^m=z^m$. Then the system of equations (\ref{eq:F-can}) takes
the form
$$\left.
\begin{array}{rrr} q^1-z^mq^{\bar{2}}-z^1 & = & 0 \\
q^1+z^mq^{\bar{1}}-z^2 & = & 0 \\ q^3-z^3 & = & 0 \\
        &\vdots &  \\ q^m-z^m & = & 0
\end{array}\right\}
$$ which has the solution
$$ z^m=q^m, \ldots , z^3=q^3,
z^2=q^2+q^mq^{\bar{1}},z^1=q^1-q^mq^{\bar{2}}.
$$ All the maps $z^k$ are globally defined and are holomorphic with
respect to the Hermitian structure $J(\mu )$.  In fact each is also
holomorphic with respect to a K\"ahler structure (and so has superminimal
fibres) and none are full if $m>3$.  However, define a holomorphic
composition of the maps $z^1, \cdots , z^m$ by
\begin{eqnarray*}
\phi & = & z^1z^2\cdots z^{m-1}\\
 & = & (q^1-q^mq^{\bar{2}})(q^2+q^mq^{\bar{1}})q^3\cdots q^{m-1}\ .
\end{eqnarray*}

By Proposition \ref{prop:sethaho} or \ref{prop:cosymplectic}, $\phi $ is
holomorphic with respect to
$J(\mu )$ and harmonic.  We show that $\phi$ does not have superminimal
fibres with respect to any almost complex structure and so is not
holomorphic with respect to any K\"ahler (or parallel complex)
structure.  Furthermore
$\phi$ is full.  For ease of exposition we consider the case $m=3$.  By
Proposition
\ref{prop:suptest},
$\nabla\phi = z^2\nabla z^1+z^1\nabla z^2$ lies in a linear subspace

\begin{equation}\label{eq:ls}
\langle b,\nabla \phi \rangle^{\cc} = 0
\end{equation} if and only if
\begin{equation} z^2\left|\begin{array}{rrr} b^1-z^3b^{\bar{2}} & 0 &
-q^{\bar{2}} \\ b^2+z^3b^{\bar{1}} & -1 & q^{\bar{1}} \\ b^3 & 0 & -1
\end{array}\right| + z^1\left|\begin{array}{rrr} -1 & b^1-z^3b^{\bar{2}}
& -q^{\bar{2}} \\ 0 & b^2+z^3b^{\bar{1}} & q^{\bar{1}} \\ 0 & b^3 & -1
\end{array}\right|=0
\end{equation}

Now
$$ q^{\bar{1}}=\frac{z^{\bar{1}}+z^2z^{\bar{3}}}{1+|z^3|^2},
q^{\bar{2}}=\frac{z^{\bar{2}}-z^1z^{\bar{3}}}{1+|z^3|^2}.
$$ Then setting $z^1z^2=a$ where $a$ is a constant and eliminating $z^1$
gives that (\ref{eq:ls}) is satisfied if and only if

\begin{eqnarray*} |z^2|^2\left[ (b^1-z^3b^{\bar{2}})(1+|z^3|^2)z^2 -
b^3(|z^2|^2-az^{\bar{3}})\right]
 & + & \\  |z^2|^2\left[ (b^2+z^3b^{\bar{1}})(1+|z^3|^2)z^{\bar{2}} +
b^3(\bar{a}+|z^2|^2+ |z^2|^2z^{\bar{3}})\right]  & = & 0.
\end{eqnarray*} Comparing coefficients, provided $a$ is non-zero, yields
$b^1=b^2=b^3=b^{\bar{1}} =b^{\bar{2}}=0$, leaving $b^{\bar{3}}$
arbitrary.  Thus along each fibre over
$a\neq 0$, the smallest subspace of $\CC^6$ containing $\nabla \phi$ is a
${\CC }^5$, showing that such a fibre is not superminimal with respect to
any almost complex structure.

To show $\phi$ is full, by Proposition \ref{prop:fulltest} , we follow
the above calculation, but now suppose $b^{\bar{k}}=\bar{b^k}$ and then
$b\equiv 0$.  For general
$m$, the smallest subspace of $\Cmm$ containing $\nabla\phi$ is a
$\CC^{m+2}$.  Thus $\phi $ provides an example of a harmonic morphism
defined globally on ${\Rmm}$ which is both full and does not have
superminimal fibres with respect to any almost complex structure.
\end{example}
  By modifying
$\phi$ to be
$$
\tilde{\phi }=z^1\cdots z^{m-1}+z^{m},
$$ we obtain a {\em submersive} example with the same properties as
$\phi$ (full, with non-superminimal fibres).  In particular the fibres of
$\tilde{\phi }$ form a global conformal foliation of ${\Rmm}$ by minimal
submanifolds of codimension 2, which does not arise from any K\"ahler
structure.  That is the leaves are not the fibres of a map, holomorphic
with respect to a K\"ahler structure.

\bigskip {\em Global examples in odd dimensions}

To construct examples on ${\RR}^{2m-1}$ we proceed as follows.

\begin{example}\label{ex:glob2} Let $\mu _i, h^i$ be defined as in
Example \ref{ex:glob1}. Define
\begin{eqnarray*}
\phi & = & (z^1z^m-z^2)z^3\cdots z^{m-1} \\
 & = & \big( (q^1-q^{\bar{1}})q^m-(q^m)^2q^{\bar{2}}-q^2\big) q^3\cdots
q^{m-1}.
\end{eqnarray*}

Then $\phi$ is invariant in the $x^1$-direction ($x^1=\hbox{Re}\, q^1$).
Calculations similar to the last example above show that along almost
every fibre, the smallest subspace of $\Cmm$ containing $\nabla\phi$ is a
$\CC^{m+1}$ and so the fibres are not superminimal with respect to any
almost complex structure.  Furthermore
$\phi$ is not invariant under any other direction apart from $x^1$ and so
reduces to a globally defined, full harmonic morphism on
${\RR}^{2m-1}$.

The modified example
$$
\tilde{\phi } = (z^1z^m - z^2) + (z^3)^2 + \cdots + (z^{m-2})^2 + z^{m-1}
$$ gives a {\em submersive} harmonic morphism on ${\Rmm}$, with fibres
not superminimal with respect to any almost complex structure, which
reduces to ${\RR}^{2m-1}$ and is full on that space.  Thus the fibres of
$\tilde{\phi }$ form a global conformal foliation by submanifolds of
codimension 2 of ${\RR}^{2m-1}$ which does not arise from any K\"ahler
structure.
\end{example}

\subsection{Results} \label{subsec:Results} From the above examples we
deduce the existence of various sorts of harmonic morphism. First the
local case:

\begin{theorem}\label{th:local-even} For any $m>2$ there are full harmonic
morphisms $\phi:U \ra \CC$ defined on open subsets $U$ of $\Rmm$ which
are holomorphic with respect to a Hermitian structure $J$ on $U$ but not
holomorphic with respect to any K\"ahler structure on $U$.  There are
such examples with
\begin{itemize}
\item the map submersive and fibres superminimal with respect to $J$ or
not superminimal with respect to any Hermitian structure,
\item $J$ cosymplectic or not,
\item holomorphic with respect to a family of Hermitian structures, with
fibres superminimal with respect to some of the family but not all.
\end{itemize}
\end{theorem}

\begin{remarks} This is not true for $m=2$ where a submersive harmonic
morphism is holomorphic with respect to precisely one Hermitian
structure, or precisely two of different orientations if the fibres are
totally geodesic, and the fibres are always superminal with respect to
those Hermitian structures.
\end{remarks}

Next the global results.  First for even dimensions.

\begin{theorem}\label{th:global-even} For any $m>2$ there are full
harmonic morphisms $\Rmm \ra \CC$ defined on  the whole of $\Rmm$ which
are holomorphic with respect to a Hermitian structure $J$ on $\Rmm$ but
not holomorphic with respect to any K\"ahler structure on $\Rmm$. These
can be chosen to be submersive.
\end{theorem}

\begin{remarks} Again this is not true for $m=2$, see \cite{Wo-3} where
it is proved that any globally defined submersive harmonic morphism
$\RR^4 \ra \CC$ is holomorphic with respect to a K\"ahler structure; it
is conjectured that this is true for any harmonic morphism whether
submersive or not.
\end{remarks}

Since the fibres of a harmonic morphism to a surface give a conformal
foliation by minimal submanifolds \cite{Ba-Ee}, we deduce
\begin{corollary}\label{cor:conf-even} For any $m>2$ there are full
conformal foliations by minimal submanifolds of codimension $2$ of the
whole of $\Rmm$ which are holomorphic with respect to a Hermitian
structure on $\Rmm$ but not with respect to any K\"ahler structure.
\end{corollary}

 Next for odd dimensions:

\begin{theorem}\label{th:global-odd} For each $m = 1,2,\ldots$ there are
harmonic maps $\Rmm
\ra \CC$ defined on the whole of $\Rmm$ which are holomorphic with
respect to a Hermitian structure $J$ on $\Rmm$ but not holomorphic with
respect to any K\"ahler structure on $\Rmm$ and which factor to full
harmonic morphisms from
$\RR^{2m-1}$ to $\CC$. These harmonic morphisms can be chosen to be
submersions.
\end{theorem}

We can combine these results for $n$ even and odd into a single statement
by using the notion of {\em arising from a K\"ahler structure}
(Definition \ref{def:arise}):

\begin{corollary}\label{cor:global-all} For any $n > 4$ there are full
harmonic morphisms from $\RR^n$ to $\CC$ which do not arise from K\"ahler
structures.
\end{corollary}

The consequence of this for foliations is

\begin{corollary}\label{cor:conf-all} For any $n > 4$ there are full
conformal foliations of $\RR^n$ by minimal submanifolds of codimension
$2$ which do not arise from K\"ahler structures.
\end{corollary}

\subsection{Harmonic morphisms from tori} A map $\Phi:\Rmm \ra \CC$
descends to a map $\phi:\Rmm/\Gamma \ra
\CC/\Gamma^{\prime}$ of tori (where $\Gamma$ and  $\Gamma^{\prime}$ are
lattices) if and only if
$\Phi(\Gamma) \subset \Gamma^{\prime}$. Then $\phi$ is a harmonic
morphism if and only if
$\Phi$ is and all harmonic morphisms $\phi$ arise in this way. The
following Proposition describes such maps.
\begin{proposition}  Let $\phi:\Rmm/\Gamma \ra \CC/\Gamma^{\prime}$ be a
harmonic morphism of tori.  Then

(i)  With respect to suitable Euclidean coordinates $q^I$ on $\Rmm$,
$\Phi = c q^1$ for some constant $c > 0$.

(ii)  For $m > 2$,  $\phi$ is holomorphic with respect to infinitely many
K\"ahler and infinitely many non-K\"ahler Hermitian structures $J$ on the
torus.  If $m > 3$, we can choose $J$ such that the fibres of $\phi$ are
not superminimal with respect to it.
\end{proposition}

\proof (i)  Since $\phi$ is harmonic, $\Phi$ must be linear, say, $\Phi =
\sum_I a_I q^I$; then, if $\phi$ is a harmonic morphism, horizontal weak
conformality (\ref{eq:HWC}) implies that the vector with components $a_I$
is isotropic.  So we can choose an orthonormal change of coordinates to
make it $(1,0,\ldots,0)$.

(ii)  The K\"ahler assertion is obvious.  For the other assertion define
$M:\Rmm /\Gamma \ra so(m,\CC)$ by
$M^2_{\bar{3}} = -M^3_{\bar{2}} =
\lambda {\cal P}(q^1)$ where ${\cal P}$ is a holomorphic function
$\CC/\Gamma^{\prime} \ra \CC
\cup \{\infty\} = \CP$ (e.g. the Weierstrass ${\cal P}$ function) and
other entries of $M$ zero.  Then set
$J = J(M)$ and note that $J$ is well-defined at poles of ${\cal P}$ since
it is of the form
$$
\Rmm / \Gamma \stackrel{q^1}{\ra} \CC / \Gamma^{\prime} \stackrel{\cal
P}{\ra} \CC
\cup \{\infty\} \stackrel{\times \lambda}{\ra} \CC \cup \{\infty\}
\stackrel{stereo}{\ra} \CP \hookrightarrow {\JJ}^{+}(\Rmm)\ .
$$ Then $J$ is easily seen to be Hermitian, and non-K\"ahler if $\lambda
\neq 0$. Note that we can see
$\phi = q^1$ as the first component $z^1$ of the solution $z(q)$  of the
equation (\ref{eq:F}) with
$M^2_{\bar{3}} = -M^3_{\bar{2}} = \lambda {\cal P}(z^1)$ and other
entries of
$M$ zero, and
$f^1 = {\cal P}(w^1) - z^1$,\ $f^i = w^i - z^i$ for $i > 1$.  For
$m>3$, we may replace $q^1$ by $q^2$ in the definition of $M$.
\eproof

{\small e-mail addresses of authors: baird@kelenn-gw.univ-brest.fr\ ;\
j.c.wood@leeds.ac.uk}


\begin{thebibliography}{99}

\bibitem{Ba-6} P. Baird, {\em Riemannian twistors and Hermitian
structures on low-dimensional space forms}, J. Math. Phys. {\bf 33}
(1992), 3340-3355.

\bibitem{Ba-Ee} P. Baird and J. Eells, {\em A Conservation law for
harmonic maps}, in: Geometry Symposium Utrecht 1980, Lecture Notes in
Math. {\bf 894}, Springer, Berlin, Heidelberg, New York, London, Paris,
Tokyo (1981), pp. 1-25.

\bibitem{Ba-Wo-6} P. Baird and J.C. Wood, {\em Harmonic morphisms and
minimal submanifolds of spheres and projective spaces}, preprint in
preparation.

\bibitem{Besse-Einstein} A.L. Besse, {\em Einstein manifolds}, Ergeb.
Math. (3) {\bf 10}, Springer, Berlin, Heidelberg, New York, London,
Paris, Tokyo, 1987.

\bibitem{BdeB} D. Burns and P. de~Bartolomeis, {\em Applications
harmoniques stables dans ${\bf P}^n$}, Ann. Sci. \'Ecole Norm. Sup.{\bf
22} (1988),159-177.

\bibitem{Ee-Sa} J. Eells and J.H. Sampson, {\em Harmonic mappings of
Riemannian manifolds}, Amer. J. Math. {\bf 86} (1964), 109-160.

\bibitem{Fu-1} B. Fuglede, {\em Harmonic Morphisms between Riemannian
manifolds}, Ann. Inst. Fourier {\bf 28} (1978), 107-144.

\bibitem{Gu-3} S. Gudmundsson, {\em Harmonic morphisms from complex
projective spaces}, Geom. Dedicata (to appear).

\bibitem{Gu-4} S. Gudmundsson, {\em Non-holomorphic harmonic morphisms
from K\" ahler manifolds}, Man. Math. (to appear).

\bibitem{Gu-Wo-1} S. Gudmundsson and J.C. Wood, {\em Multivalued Harmonic
Morphisms}, Math. Scand. {\bf 73} (1994), 127-155.

\bibitem{Gu-Wo-2} S. Gudmundsson and J.C. Wood, preprint in preparation.

\bibitem{Is} T. Ishihara, {\em A mapping of Riemannian manifolds which
preserves harmonic functions}, J.Math.Kyoto Univ. {\bf 19} (1979),
215-229.

\bibitem{Jost-Yau} J. Jost and S.-T. Yau, {\em A nonlinear elliptic
system for maps from Hermitian to Riemannian manifolds and rigidity
theorems in Hermitian geometry}.

\bibitem{Lichn} A. Lichnerowicz, {\em Applications harmoniques et
vari\'et\'es k\"ahl\'eriennes}, Rend. Sem. Mat. Fis. Milano {\bf 39}
(1969), 186-195.

\bibitem{Metzler-Muir} W. Metzler and T.H. Muir, {\em A treatise on the
theory of determinants}, Dover Publications Inc., New York, 1933.

\bibitem{New-Nir} A. Newlander and L. Nirenberg, {\em Complex analytic
coordinates in almost complex manifolds}, Ann. of Math. {\bf 65} (1957),
391-404.

\bibitem{Wo-3} J.C. Wood, {\em Harmonic morphisms and Hermitian
structures on Einstein 4-mani\-folds}, Internat. J. Math. {\bf 3} (1992),
415-439.

\end{thebibliography}
\end{document}